\documentclass{article}
\usepackage{amsmath}
 \usepackage{amscd}
 \usepackage{amssymb}
\usepackage{graphicx}
 \usepackage{mathptmx}      
\usepackage{amsmath}
 \usepackage{amscd}
 \usepackage{amssymb}
\usepackage{amssymb}
  \usepackage[all]{xy}

\newtheorem{TEO}{Theorem}[section] 
\newtheorem{PROP}[TEO]{Proposition}
\newtheorem{LEM}[TEO]{Lemma}
\newtheorem{DEF}[TEO]{Definition}
\newtheorem{COR}[TEO]{Corollary}

\newtheorem{defa}{Definition}

\newcommand\Oh{{\mathcal O}}
\newcommand\sA{{\mathcal A}}

\newcommand\sF{{\mathcal F}}
\newcommand\sG{{\mathcal G}}
\newcommand\sH{{\mathcal H}}
\newcommand\sI{{\mathcal I}}

\newcommand\sK{{\mathcal K}}

\newcommand\sQ{{\mathcal Q}}
\newcommand\sR{{\mathcal R}}

\newcommand\om{\omega}

\newcommand\Ga{\Gamma}

\newcommand\simby[1]{\buildrel{{\scriptstyle\mathrm{#1}}}\over{\sim}}
\newcommand\numeq{\simby{{num}}}         
\newcommand\dual{\mathrel{\raise3pt\hbox{$\underline{\mathrm{\thinspace d
\thinspace}}$}}}

\newcommand\iso{\cong}

\newcommand\ra{\rightarrow}
\newcommand\into{\hookrightarrow}
\newcommand\onto{\twoheadrightarrow}

\newcommand\proj{\mathbb P}
\newcommand\Ka{\mathbb K}
\newcommand\Na{\mathbb N}
\newcommand\Z{\mathbb Z}

\newcommand\coker{\operatorname{coker}}
\newcommand\im{\operatorname{im}}

\newcommand\mult{\operatorname{mult}}

\newcommand\Pic{\operatorname{Pic}}

\begin{document}

	
		\title{On the canonical ring of \\
curves and surfaces 
\thanks{Research carried out under the
MIUR  project 2008  ``Geometria delle Variet\`a algebriche e dei loro spazi dei moduli''}}
	 \author{Marco Franciosi}

	\date{}

\maketitle

	\begin{abstract} 
 Let $C$ be  a 
 curve   (possibly non reduced or reducible) 
 lying on a smooth algebraic surface.
 We show that the canonical ring 
 $  R(C, \omega_C)=\bigoplus_{k\geq 0} H^0(C, {\omega_C}^{\otimes k})$
  is generated in degree 1 if $C$ is numerically  $4$-connected, not 
  hyperelliptic and even (i.e. with $\omega_C$ of even degree on every component). 
  
   As a corollary we show that on a smooth algebraic surface of general type 
   with $p_g(S)\geq 1$ and $q(S)=0$
   the
   canonical ring  $  R(S, K_S)$ is generated in degree $\leq 3$  if 
  there exists a curve $C$  $ \in |K_S|$  numerically $3$-connected and not 
  hyperelliptic.

 \hfill\break	
{\bf keyword: }
algebraic curve,  surface of general type, canonical ring,  pluricanonical embedding.  

\hfill\break  {\bf Mathematics Subject Classification (2000)} 14H45,  14C20,  14J29
\end{abstract}

\section{Introduction}
 \label{intro}

Let  $C$ be a curve  
 (possibly non reduced or reducible) 
 lying on a smooth algebraic surface $S$   and let 
  $\om_C$ be the   dualizing  sheaf of $C$. 
 The purpose of this paper is to analyze the  canonical ring of $C$, that is,  the 
 graded ring
$$ \displaystyle{ R(C, \omega_C)=\bigoplus_{k\geq 0} H^0(C, {\omega_C}^{\otimes k})}$$ 
under some suitable assumptions on the curve $C$. 

The rationale  of our  analysis stems from  
 several aspects of  the theory of algebraic surfaces.

The first  such aspect is the analysis  of  surface's fibrations and the study of
their applications to  surface's geography.  Indeed, 
given a  genus $g$ fibration $ f: S\rightarrow B$  over a smooth curve $B$, 
an important tool in this analysis is the  {\em  relative 
canonical algebra}   $ R (f)= \bigoplus_{n\geq 0} f_{\ast}(\omega_{S/B}^{\otimes n})$. 

In recent years  the importance of  $ R (f)$ has become clear (see  Reid's  unpublished note  \cite{Reid}) 
and a way to understand its behavior consists in  studying  the canonical ring of 
every fibre of $f$.   
More specifically, 
   denoting
  by  $ C = f^{-1}(P)$  the scheme theoretic fibre
over a point $p\in B$,
the local structure around  $P$ of  the {\em  relative 
canonical algebra}  
can be understood via the canonical algebra of $C$, since     
  the reduction modulo  ${\cal M}_P$  of 
the stalk at $P$ of the {\em  relative 
canonical algebra}   is  nothing but   $R(C, \omega_C)$  (see  \cite[\S1]{Reid}).
 Mendes Lopes in \cite{ML}  studied the cases   where the genus  $g$ of the fibre is $g\leq 3$ whereas 
in \cite{Kon} and \cite{Fr} it is shown that  for every $g \geq 3$,  
$R(C, \omega_C)$ is generated in degree $\leq 4$
 if  every fibre is numerically connected 
 and  in degree  $\leq 3$  if furthermore  there are no multiple fibres.

More recently Catanese and Pignatelli in \cite{CP}  illustrated 
 two  structure theorems for fibration of low genus
using a detailed description 
of the relative canonical algebra. In particular 
they  
showed  an interesting characterization of isomorphism classes of 
relatively minimal fibration of genus $2$  and of relatively minimal fibrations of genus $3$
with fibres numerically 2-connected  and not hyperelliptic
 (see \cite[Thms.  4.13, 7.13]{CP}).

 Another motivation  of  our analysis  (see \cite{Lau2}) relies  on the study  of 
the resolution of a normal  surface singularity $\pi:S\rightarrow X$. Indeed 
 if one considers the fundamental cycle arising from a minimal resolution 
	 of a singularity,  then  the study of the ring
	 $ \displaystyle{ R(C, K_S)=\bigoplus_{k\geq 0} H^0(C, {K_S}^{\otimes k})}$ plays a  
	 crucial role in order to understand that singularity.

Finally,
 as shown in \cite{CFHR},   the study of invertible sheaves on curves  possibly reducible or non reduced
is rich in implications 
 in the cases  where Bertini's  theorem  does not hold or  simply  
if one needs to consider every curve contained  in a given linear system.  
 For instance,   one can acquire  information on the canonical ring of a surface of general type 
 simply by  taking its 
 restriction
to an effective canonical  divisor $C \in |  K_S|$  (not necessarily irreducible, neither 
reduced)  and considering the canonical ring $R(C, \omega_C)$ (see Thm. \ref{teo:bicanonical} below).

\hfill\break 
In this paper we  
 analyze the  canonical ring of $C$ 
when the curve $C$ is  $m$-connected and {\em even}, and we  show some applications to 
the study of the canonical ring of an algebraic surface of general type.

For a curve $C$ lying  on a smooth algebraic surface $S$, being {\em  m-connected}  means that 
$C_1 \cdot C_2 \geq m $    
for every effective decomposition  $ C=C_1 + C_2$,
(where $C_{1} \cdot C_{2}$  denotes their intersection number as 
divisor on $S$). 
If $C$ is $1$-connected usually $C$ is said to be numerically connected.  
The definition turns back to Franchetta (cf. \cite{frc})
and has many relevant implications. For instance  in \cite[\S 3]{CFHR} it is shown that if the curve $C$ is 
$1$-connected then $h^0( C, \Oh_C)=1$, if $C$ is $2$-connected then the system $|\omega_C|$ is base point free,
whereas if $C$ is $3$-connected  and not honestly hyperelliptic 
(i.e.,   a finite double cover of 
 $\proj^1$ induced by the canonical morphism) 
 then $\omega_C$ is very ample. 

Keeping the usual notation for effective divisor 
 on  smooth  surfaces, i.e.,    writing $C$  as
  $\sum_{i=1}^{s} n_{i} \Gamma_{i}$ 
(where  
 $\Gamma_{i}$'s     are  the 
 irreducible components of $C$  and  $n_{i}$   the multiplicity of 
 $\Gamma_{i}$ in $C$),  
the second condition  can be illustrated by the following definition.

\begin{defa}\label{def:evencurve} 	Let $C=\sum_{i=1}^{s} n_{i} \Gamma_{i}$ 
	be a  curve 
	contained in a smooth algebraic surface.   $C$ is said to be 
	{\em even} if 
	$   \deg ({\om_{C}}_{| \Gamma_{i}}) $ is even for every irreducible 
	$ \Gamma_{i} \subset C $
	(that is,   $\Gamma_{i}\cdot (C -\Gamma_{i})$ even for every
	$ i=1,...,s$.) 
\end{defa}
Even curves  appear for instance when  considering  the 
canonical system $|K_{S}|$ for a surface $S$ of general type.  Indeed, 
by adjunction, for every curve $C \in |K_{S}|$ we have 
$|(2K_{S})_{|C}|=
|K_{C}|$, that is,  every curve in the 
canonical system is even.

The main result of this paper is a generalization to even curves of  the classical 
Theorem of Noether and 
Enriques on the degree of the generators of the graded ring $R(C, \omega_C)$: 
\begin{TEO}\label{teo:3}
 Let  $C$   be an even    $4$-connected  curve 
contained in a smooth algebraic surface. If  $p_{a}(C) \geq 3$
and  $C$ is  not honestly hyperelliptic 
then $R(C, \omega_C)$ is generated in degree 1.
\end{TEO}
Following the notations of  \cite{Mu}, this result can be rephrased by saying that 
$\omega_{C}$ is normally generated on $C$.  In  this case  the embedded curve 
$\varphi_{|\omega_C|}( C) \subset \proj^{p_a(C) -1}$ is 
arithmetically   Cohen--Macaulay. 
 
The proof of Theorem \ref{teo:3}  is based on  the ideas 
adopted by Mumford in  
  \cite{Mu} and on the results obtained in \cite{Fr} for 
adjoint divisors, via a detailed analysis of the possible decompositions of the given curve $C$.

As a corollary we obtain a bound on the degree of the generators of  the   canonical ring of a surface of general type.

If  $S$ is a smooth  algebraic surface and $K_S$ a canonical divisor, the canonical ring of $S$ is
the graded algebra 
 $$ R(S, K_S)=\bigoplus_{k\geq 0} H^0(S, {K_S}^{\otimes k})$$ 
In \cite{ci} a detailed analysis of $ R(S, K_S)$ is   presented   in the most interesting case where
$S$ is of general type and  there are given bounds (depending
on the invariants $p_g(S) := h^0(S, K_S)$, $q:= h^1(S, \Oh_S)$,  and $K_S^2$)   on the degree 
of elements of $R(S, K_S)$  forming a minimal system of homogeneous generators.  
Furthermore it is shown that for small values of $p_g$ some exceptions  do occur, 
depending substantially  on the numerical connectedness of the curves in the linear system $|K_S|$. 
In particular  in \cite[\S 4]{ci} there are given examples of surfaces of general type   with
$K_S$ not  $3$-connected  whose  canonical ring is not generated in degree $\leq 3$ and it is 
conjectured that the $3$-connectedness of the canonical divisor $K_S$ should imply
 the generation of $R(S, K_S)$
 in degree
1,2,3, at least in the case $q=0$.

Here we show that  this is the case. Our result, obtained essentially  by restriction to a curve
$C\in |K_S|$, is the following
\begin{TEO}\label{teo:bicanonical}
	Let $S$ be a surface of general type 
	with  $p_g(S) := h^0(S, K_S) \geq 1$ and  $q:= h^1(S, \Oh_S)=0$. Assume that there exists a curve 
	$C\in
|K_S|$ such that $C$ is numerically $3$-connected and not 
 honestly hyperelliptic. Then  
 the canonical ring of $S$ 
 is generated in degree $\leq 3$. 

 \end{TEO}

The paper is organized as follows: in 
\S  2 some useful background results are illustrated; 
 in \S 3 we point out a result on 
 even divisors; 
 in \S 4 we consider the special case of binary curves, i.e. the case where $C$ is the union of two rational curves;
   in \S 5 we prove Thm. \ref{teo:3}; in \S 6 we give the proof of Thm. \ref{teo:bicanonical}.

\hfill\break
{ \em Acknowledgment. }  The author  wish to thank  the Mathematisches Institut 
 of the
 University of Bayreuth  (Germany) for its support and its warm  hospitality during the preparation of part of this work.

The author would like to thank Elisa Tenni for her careful reading and the referee  for its  stimulating suggestions.

\section{Notation and preliminary results}
 \subsection{Notation} 
 
 We work over  an algebraically closed field $\Ka$ of characteristic $\geq 0$.
 
Throughout this paper $S$ will be a smooth algebraic surface         
over $\Ka$ and  $C$ will be a 
curve lying on $S$ (possibly reducible and non
reduced). Therefore $C$ 
will be written 
(as a divisor on $S$) as $C=\sum_{i=1}^{s} n_{i} \Gamma_{i}$, where 
the  $\Gamma_{i}$'s are    the 
 irreducible component of $C$  and the $n_{i}$'s are   the multiplicities. 
A subcurve $B\subseteq C$  
will mean a curve 
 $\sum m_{i} \Gamma_{i}$, with $0\leq m_i\leq n_i$ for every $i$.

By abuse of notation if $B \subset C $ is a subcurve  of $ C$,   $ C-B$  denotes the  curve $A$ such that 
 $C= A+B$ as divisors on $S$.

\hfill\break 
 Let $\sF$ be an invertible sheaf on $C$.  

 If $G\subset C$ is a proper subcurve of $C$ then ${\sF}_{| G} $ denotes its restriction to $G$. 
 
 For each $i$ the natural inclusion map $\epsilon_i:\Ga_i \rightarrow C$ 
induces a map 
$\epsilon_i^{\ast}:\sF \rightarrow \sF_{|\Ga_i}$. We denote by  
 $d_i= \deg (\sF_{|\Ga_i}) = \deg_{\Gamma_i} \sF$ the degree of $\sF$ on each  irreducible component,
 and by  ${\bf d}:=(d_1, ..., d_s)$  the 
 multidegree  of $\sF$ on $C$.  If  $B = \sum m_{i} \Gamma_{i}$ is a subcurve of $C$, by $\mathbf 
d_{B}$  we mean the multidegree of $\sF_{|B}$.

  By $\Pic^{\mathbf d}(C)$ we  denote the Picard scheme which 
parametrizes the classes of invertible
sheaves of multidegree ${\bf d}=(d_{1},\ldots,d_{s})$ (see \cite{Fr}).

We recall that for every ${\bf d}=(d_1, ..., d_s)$ there is an isomorphism
$\Pic^{\bf d}(C)\iso \Pic^{\bf 0}(C)$ and  furthermore 
$\dim \Pic^{\bf 0}(C) = h^1(C, \Oh_C)$  (cf. e.g.  \cite{BPV}).  

Concerning the Picard group of $C$ and the Picard group
 of a subcurve $B \subset C$
we have  
$$\Pic^{\mathbf d}(C) \onto 
\Pic^{\mathbf d_{B}}(B)  \  \  \  \forall {\mathbf d}$$ 
(see  \cite[Rem. 2.1] {Fr}). 
 
An invertible sheaf  $\sF$ is  said to be NEF if $d_i \geq 0$ for every $i$.  
Two invertible sheaves $\sF$, $\sF'$ are said to be 
 numerically
 equivalent on $C$ 
 (notation: $\sF\numeq \sF'$)  if 
 $\deg_{\Gamma_{i}} \sF=\deg_{\Gamma_{i}}  {\sF'} $
 for every   $ \Gamma_{i} \subseteq  C$. 

\hfill\break 
$\om_C$ denotes the   dualizing sheaf of $C$ (see \cite{Ha}, Chap.~III, \S7), and 
$p_a(C)$ the arithmetic genus of $C$, $p_a(C)=1-\chi(\Oh_C)$.

 If $G\subset C$ is a proper subcurve of $C$  we denote by $H^0(G, \omega_C)$
 the space of sections of  $ {\omega_C}_{|G}$.

\hfill\break 
Finally,   a  curve 
 $C$ is  said to be {\em honestly
hyperelliptic}  if there exists a finite
morphism $\psi\colon C\to\proj^1$ of degree $2$. In this case 
 $C$ is either irreducible, or of the form $C=\Gamma_{1}+
\Gamma_2$ with $p_{a}(\Gamma_{i} ) = 0$ and 
$\Gamma_1 \cdot \Gamma_2=p_{a}(C) +1$
(see    \cite[\S3]{CFHR} for a detailed treatment).

 \subsection{General divisors of low degree}

Let 
$C=\sum_{i=1}^{s} n_{i} \Gamma_{i}$ be a curve 
 lying on a smooth algebraic surface $S$.   
  An invertible sheaf  on $C$ of multidegree  ${\bf d} =(d_1, ..., d_s) $
 is said to be 
 ``general''  if the corresponding  class in the Picard 
scheme  $\Pic^{\mathbf d}(C)$ is in general position, i.e.,  if it lies
in the complementary of a proper closed subscheme (see \cite{Fr} for details).  

We recall two vanishing  results 
 for general invertible sheaves  of low degree.

 \begin{TEO}[  {
 \cite[Thms. 3.1, 3.2 ]{Fr} }]\label{teo:general} 
 
  \begin{enumerate} 
 \renewcommand\labelenumi{(\roman{enumi})}
\item  If $\sF$ is 
	     a ``general''  invertible  
 sheaf    such that  
 $ \deg_{B} \sF \geq \ p_a(B)$  for every  subcurve $B\subseteq C$, then 
 it is $H^{1}(C,\sF)=0$.  
\item 	If $\sF$  is 
a ``general'' invertible sheaf 
such that 
$  \deg_{B} \sF \geq \ p_a(B) + 1$   for every  subcurve $B\subseteq 
C$, 
 then  the linear system $|\sF|$ 
is base point free.
\end{enumerate}
\end{TEO} 
In particular   we obtain the following 

\begin{PROP}\label{mezzocanonico} Let $C = \sum_{i=1}^{s} n_{i}\Gamma_{i}$  be  
a  $1$-connected curve contained in a smooth algebraic surface, and consider a proper subcurve $B\subsetneq C$.
Let 
${\bf d}= (d_1, ..., d_s) \in \Z^{s} $ be such that 
$d_i \geq \frac{1}{2} \deg_{\Gamma_i} \omega_C $
$\forall \ i =1,...,s$. 

Then for  a ``general"  invertible sheaf  $\sF$ in $ \Pic^{\mathbf d_{B}}(B)$
 it is:
\begin{enumerate}
	 \renewcommand\labelenumi{(\roman{enumi})}
\item $H^1(B, \sF) = 0 $; 

\item  $| \sF_{| B}| $ is a base point free system on $B$ if $C$ is 3- connected.

\end{enumerate}
\end{PROP} 
Considering the case where $C$ is  an  even  $4$-connected  curve   we obtain 
\begin{COR}\label{mezzocanonico2} Let $C = \sum_{i=1}^{s} n_{i}\Gamma_{i}$  be  
a  $4$-connected  even curve contained in a smooth algebraic surface. 

For every $i=1,\dots,s$,  let 
$d_i=\frac{1}{2}  \deg_{\Gamma_i} \omega_C$  and let  ${\bf d}= (d_1, ..., d_s) \in \Z^{s} $.

Let $B\subsetneq C$  be  a proper subcurve  of $C$  and consider a   
 a ``general"  invertible sheaf $\sF$   in $ \Pic^{\mathbf d_{B}}(B)$ (i.e.,  with an abuse of notation  we can write
 $\sF \numeq \frac{1}{2} {\omega_C}_{|_{B}}$).

Then 
$H^1(B, \sF) = 0 $ and  $|\sF_{| B}|$ is a base point free system.

\end{COR}

\subsection{Koszul cohomology groups of algebraic curves} 
Let $C = \sum_{i=1}^{s} n_{i}\Gamma_{i}$ be a curve 
 lying on a smooth algebraic surface $S$ and   let 
$\sH, \ \sF$ be  invertible sheaves on $C$. Consider   a subspace $W \subseteq 
H^{0}(C,\sF)$  which yields a base point free system of 
projective dimension $r.$

The Koszul groups $\sK_{p,q}(C, W,\sH, \sF)$ are defined 
as the cohomology at the middle of the complex 
\[
\stackrel{p+1}{\bigwedge} W \otimes H^{0}(\sH\otimes \sF^{q-1})  \longrightarrow 
\stackrel{p}{\bigwedge} W \otimes H^{0}(\sH\otimes \sF^{q})  \longrightarrow 
\stackrel{p-1}{\bigwedge} W \otimes H^{0}(\sH\otimes \sF^{q+1})  
 \] 
If $W=H^{0}(C,\sF) $ they are usually denoted by 
$\sK_{p,q}(C, \sH, \sF)$, while if $\sH \iso \Oh_C$ the usual notation is $\sK_{p,q}(C, \sF)$ 
(see \cite{Gr} for  the  definition and main results).

We point out  that the multiplication map 
\begin{equation}\nonumber
 W \otimes H^0(C,\sH) \rightarrow H^0(C,\sF\otimes \sH) 
 \end{equation} 
is surjective
if and only if
 $\sK_{0,1}(C, W, \sH, \sF) =0$ and the ring 
  $R(C, \sF)=\bigoplus_{k\geq 0} H^0(C, {\sF}^{\otimes k})$ is generated in degree $1$
   if and only if $\sK_{0,q}(C, \sF) =0$ $\forall \ q \geq 1 $. 
Moreover  if $\sF$ is very ample and  $R(C, \sF)$ is generated in degree $1$, then, identifying 
$C$ with its image in $\proj^r \iso \proj(H^0(\sF)^{\vee})$, 
it is  $ \sK_{1,1}(C, \sF) \iso I_2(C, \proj^r )$, the space of quadrics 
in $\proj^r$ vanishing on $C$ (see \cite{Gr}). 

For our analysis 
the main applications of Koszul cohomology 
 are    
 the following propositions (see \cite[\S 1]{Fr}, \cite[\S 1]{Kon} for further details on curves lying on smooth surfaces).

\begin{PROP}[Duality]\label{duality}  Let $\sF$, $\sH$ be  invertible sheaves 
on $C$ and assume   $W \subseteq H^{0}(C, \sF) $ 
to	be a subspace of $\dim = r+1$ which yields 
	a base point free system. Then 
$$\sK_{p,q}( C, W, \sH, \sF) \dual 
\sK_{r-p-1,2-q}(C, W, \omega_{C}\otimes \sH^{-1}, \sF) $$
(where $\dual$ means duality of vector space).
\end{PROP}
For a proof see \cite[Prop. 1.4]{Fr}.   Following 
  the ideas outlined in   \cite[Lemma 1.2.2]{Kon} we have  a
slight generalization of Green's $H^0$-Lemma. 
 \begin{PROP}[$H^{0}$-Lemma]\label{prop:van} Let $C$ be 1-connected and 
let $\sF$, $\sH$ be  invertible sheaves 
on $C$ and assume   $W \subseteq H^{0}(C, \sF) $ 
to	be a subspace of $\dim = r+1$ which yields 
	a base point free system.
If either 
 \begin{enumerate}
	 \renewcommand\labelenumi{(\roman{enumi})}
	 \item 
$H^1(C,  \sH\otimes \sF^{-1})=0$, \\ or   
\item
$C$ is numerically connected,
$ \om_{C} \iso \sH\otimes \sF^{-1}$  and $r \geq 2$,\\  or
\item 
$C$ is numerically connected,
$h^{0}(C, \omega_{C}\otimes \sH^{-1}
 \otimes \sF ) \leq r -1$ and 
there exists a  reduced subcurve $B\subseteq C$ such that:
\begin{enumerate}
\item[$\bullet$] $ W \iso W_{|B}  $, 
\item[$\bullet$]  $ H^{0}(C, \omega_{C}\otimes \sH^{-1}
 \otimes \sF ) \into H^{0}(B, \omega_{C}\otimes \sH^{-1}
 \otimes \sF ) $,    
\item[$\bullet$]   every non--zero section of 
 $H^{0}(C, \omega_{C}\otimes \sH^{-1}
 \otimes \sF ) $ does not vanish identically  
 on any component of $B$;
 \end{enumerate}
 \end{enumerate}
then $\sK_{0,1}(C, W,  \sH , \sF)=0 $, that is, the multiplication map
$$ W \otimes H^0(C,\sH) \rightarrow H^0(C,\sF\otimes \sH) $$ is surjective.
\end{PROP}
{\bf Proof.}   
By duality we need to prove that
$
\sK_{r-1,1}(C, W, \omega_{C}\otimes \sH^{-1} , \sF)=0.$ To this aim 
let $\{ s_{0},\ldots, s_{r}\}$ be a basis for $W$ and let
 $\alpha = \sum 
 s_{i_{1}}\wedge s_{i_{2}}\wedge \ldots \wedge 
s_{i_{r-1}} \otimes \alpha_{i_{1}i_{2}\ldots i_{r-1}} \in \bigwedge^{r-1} W 
\otimes H^{0}(C, \omega_{C}\otimes \sH^{-1}
 \otimes \sF )  $
 be an element in the Kernel of the Koszul map $ 
 d_{r-1,1}$. 
 
 In cases (i)   obviously $\alpha=0$ since by Serre duality it is $H^{0}(C, \omega_{C}\otimes \sH^{-1}
 \otimes \sF ) \iso H^1(C,  \sH\otimes \sF^{-1})=0$.  In case  and (ii)  it is $H^{0}(C, \omega_{C}\otimes \sH^{-1}
 \otimes \sF )= H^0(C, \Oh_C)=\Ka$ by connectedness  and we conclude similarly
   (see also \cite[Prop. 1.5]{Fr}).
 
In the latter case by our assumptions we can restrict to the curve $B$.
Since $B$ is reduced we can choose $r+1$ ``sufficiently general points'' 
on $B$ so that $s_{j}(P_{i})=\delta^{i}_{j}$. 
But then $\alpha \in \ker  
( d_{r-1,1})$ 
  implies for every multiindex ${\bf I}=\{i_{1},\ldots i_{r-2}\}$ 
 the following equation (up to sign)
 $$  \alpha_{j_{1}i_{1}\ldots i_{r-2}} \cdot s_{{j_{1}}} + 
 \alpha_{j_{2}i_{1}\ldots i_{r-2}}\cdot s_{{j_{2}}}  +
  \alpha_{j_{3}i_{1}\ldots i_{r-2}} s_{{j_{3}}} =0. $$
 (where $\{i_{1},\ldots i_{r-2}\} \cup \{j_{1},j_{2},j_{3}\} = \{0,\ldots, 
 r+1\}$).
 
 Evaluating at $P_{j} 's$ and reindexing  we get 
 $  \alpha_{i_{1}\ldots i_{r-1}}(P_{i_{k}}) = 0 $ for  $k= 
 1,\ldots, r-1.$ 
 
 Let $\tilde{r} =  h^{0}(C, \omega_{C}\otimes \sH^{-1}
 \otimes \sF )$.  
Since the $P_{j} 's$ are in general positions and 
  every section of 
 $H^{0}(C, \omega_{C}\otimes \sH^{-1}
 \otimes \sF ) $ does not vanish identically on any component of $B$, 
 we may  assume that  any $(\tilde{r} -1)$-tuple of points $P_{i_{1}},\ldots, P_{i_{r-1}}$
 imposes 
 independent conditions   on 
 $H^{0}(C, \omega_{C}\otimes \sH^{-1}
 \otimes \sF )  $. 
 The proposition then follows by a dimension count since by 
 assumptions
 $\tilde{r} =h^{0}(C, \omega_{C}\otimes \sH^{-1}
 \otimes \sF ) \leq h^{0}(C, \sF) -2= r-1.$

\hfill$\Box$

In some particular cases we can obtain deeper  results, which 
will turn out to be useful for our induction argument in the proof of 
Theorem \ref{teo:3}. 

\begin{PROP}\label{cor:irreducible}
	Let $C$ be either 
	\begin{enumerate}
	 \renewcommand\labelenumi{(\roman{enumi})}
\item 
an irreducible curve of arithmetic genus $p_{a}(C) \geq 1$; 	 
 \\ or   
\item 
 a binary curve of genus $\geq 1$, that is,  
 $C = \Gamma_{1} + \Gamma_{2} $, with $\Gamma_{i}$ irreducible and reduced 
 rational curves 
s.t.
 $\Gamma_1 \cdot \Gamma_{2} = p_{a}(C) + 1 \geq 2	$ (see \S 4 for details).

 \end{enumerate}
Let
	$\sH \numeq \omega_{C}\otimes \sA$  be a very  
	ample divisor on $C$ s.t. $\deg_{C}  \sA \geq 4$.
	 
	Then  $\sK_{0,1}(C,  \sH , \omega_C)=0 $, that is  $H^{0}(C, \om_{C}) \otimes H^0(C,\sH) \onto 
	H^0(C,\om_{C} \otimes \sH) $. 
\end{PROP}
{\bf Proof.}
	If $p_{a}(C) =1$  then  under our assumptions it is $\om_C \iso \Oh_C$, whence the theorem follows 
	easily.  
	
	If $p_{a}(C) \geq 2$ then  by \cite[Thms. 3.3, 3.4]{CFHR} $|\om_{C}|$ is b.p.f.  and moreover it is 
	very ample if $C$ is not honestly hyperelliptic. We apply    Prop.
	\ref{prop:van}
	 with $\sF = \omega_{C}$ and $W= H^{0}(\om_{C})$.

	If $C$ is irreducible  and $ h^{0}(C, \omega_{C}\otimes \sA^{-1} )=0$ then the result follows  by  ({\em  i}) of  Prop. \ref{prop:van}. If 
	$ h^{0}(C, \omega_{C}\otimes \sA^{-1} ) \neq 0$ 
	and $h^{0}(C , \sA) = 0$	  it follows by 
		Riemann-Roch.   In the remaining case  we obtain $ h^{0}(C, \omega_{C}\otimes \sA^{-1} )\leq h^{0}(C, \om_{C}) -2$ 
	by Clifford's theorem since $\deg_{C}  \sA \geq 4$.

If $C= \Gamma_1 +\Gamma_2$ and $p_a(C) \geq 2$ we consider firstly  the case 
	where  $ \deg_{ \Gamma_{i}}  \sA_ \geq -1$  for $i=1,2$. 
	Under this assumption 
 any 
	non-zero section  of $H^{0}(C, \omega_{C}\otimes \sA^{-1} )$  does not 
	vanish identically on any single component of $C$  (otherwise 
		it would yield a section in
	$H^{0}(\Gamma_{i}, \omega_{\Gamma_{i}}\otimes \sA^{-1} ) \iso
	H^{0}(\proj^{1}, - \alpha)$ with $\alpha \geq 1$). Therefore
  we can proceed  exactly  as in the irreducible case.

Now assume   $C= \Gamma_1 +\Gamma_2$, 
$\deg_{ \Gamma_{2}}  \sA\leq  - 2$ and  $\deg_{ \Gamma_{1}}  \sA
\geq 6$.  In this case we can     apply ({\em iii}) of  Prop.
	\ref{prop:van}  taking $B=\Gamma_{2}$. Indeed,  it is 
	$h^{0}(\Gamma_{1}, \om_{\Gamma_{1}}\otimes \sA^{-1} ) =  
	h^{0}(\Gamma_{1}, \om_{\Gamma_{1}} ) =  0$ 
and 
	we have the 
	following maps 	
	$$	H^{0}(C , \omega_{C}) \iso 	H^{0}(\Gamma_{2} , \omega_{C})
	 \hspace{0,1 cm} ;  \hspace{0,2 cm}
	H^{0}(C , \omega_{C}\otimes\sA^{-1}) \into 	H^{0}(\Gamma_{2} ,
\omega_{C}\otimes\sA^{-1}).$$ 
To complete the proof it remains  to show that 
$ h^{0}(C, \omega_{C}\otimes \sA^{-1} ) \leq h^{0}(C, \omega_{C}) - 
2 = p_a(C)-2$. This follows by the following exact sequence 	
	\[ 0\rightarrow 
 H^0(\Gamma_{2}, \om_{\Gamma_{2}}\otimes \sA^{-1} )  \rightarrow
 H^0( C,  \om_{C}\otimes \sA^{-1}) \rightarrow
 H^0( \Gamma_{1} , \om_{C} \otimes \sA^{-1})  \rightarrow 0.
\]
In fact if $\deg_{\Gamma_{1}} (\om_{C} \otimes \sA^{-1}) \geq 0$ we 
have  	$ h^{0}(C, \omega_{C}\otimes \sA^{-1} ) = p_{a}(C) - 1 - \deg 
\sA$, whereas 
it is 
$ h^{0}(C, \omega_{C}\otimes \sA^{-1} ) = 
h^{0}(\Gamma_{2}, \om_{\Gamma_{2}}\otimes \sA^{-1} ) = - 
\deg_{\Gamma_{2}} \sA - 1 <  p_{a} (C) -2 $  if $\deg_{\Gamma_{1}}(\om_{C} \otimes \sA^{-1})  < 0$ 
since $ \deg_{\Gamma_{2}} (\om_C + \sA) \geq 1  $  by 
the ampleness of 
 $\om_{C} \otimes \sA$.   

\hfill$\Box$   

If one considers a curve $C$ with many components 
another useful tool is the 	following 
long exact sequences for Koszul groups.

\begin{PROP}\label{splitting} Let $C= A+B$  and  let
$|\sF|$ be a complete base point free system on $C$ such that  $H^0(C, \sF ) \onto H^0(A, \sF  )  $, 
$H^0(C, \sF^{\otimes k}  ) \onto H^0(B, \sF^{\otimes k}  )  $  
 for every  $k \geq 2$  
and 
$H^0(A, \sF  (- B)) =0$.
Then we have a long exact sequence 
$$ 
\begin{array} {llclll}
\cdots  & \to & \sK_{p+1,q-1}(C, \sF) & \to \sK_{p+1,q-1}(B,W, \sF)& \to \sK_{p,q}(A, \Oh_A(-B) , \sF) \\
&  \to & \sK_{p,q}(C, \sF) & \to \sK_{p,q}(B,W, \sF) & \to \cdots  
\end{array}
$$
where $W=H^0(C, \sF)_{|B}= \im \{H^0(C, \sF) \to H^0(B, \sF) \}$.

\end{PROP}
{\bf Proof.}
Let us consider 
$$ B^1=\bigoplus_{q\geq 0} H^0(A, \sF^{\otimes q}  (-B) ) , \  \ B^2=\bigoplus_{q\geq 0} H^0(C, \sF^{\otimes q} ) , \  \ 
B^3= W \oplus ( \bigoplus_{q\neq 1 } H^0(B, \sF^{\otimes q} )) $$
By our  hypotheses the above vector spaces can be seen as $S(H^0(C, \sF))-$modules and moreover they fit into 
the following exact sequence 
 $$ 0 \to \  B^1 \to \ B^2 \to \ B^3 \to 0  $$  
  where the  maps preserve the grading. By the long exact sequence for Koszul Cohomology
  (cf. \cite[Corollary 1.4.d, Thm. 3.b.1 ]{Gr}) we can conclude.

\hfill$\Box$ 
		  
Notice that if $C$ is numerically connected, $\sF \iso \om_C$,  $B$ is numerically connected and  $A $ is  
 the disjoint union of  irreducible rational curves 
 then the above hypotheses are satisfied. 
 
We point out   that if $\sF $ is very ample but the restriction map $H^0(C, \sF) \to H^0(B, \sF) $ is not surjective  then, 
following the notation of \cite{bf},   we can talk of ``Weak Property $N_p$" for 
 the curve $B$ embedded by the system  $W=H^0(C, \sF)_{|B}$.

\subsection{Divisors normally generated on algebraic curves} 
To conclude this preliminary section  we recall a theorem  proved in \cite{Fr} 
which allows us to start our analysis.

\begin{TEO}
[  {
 \cite[Thm. A]{Fr} }]\label{teo:adjoint} 
 Let  $C$   be a   curve 
contained in a smooth algebraic surface
and
let 
 $\sH \numeq \sF\otimes \sG$, where $\sF, \sG$ are   invertible sheaf
  such that
 $$ \begin{array}{ll}
\  \deg \sF_{|B} \geq \ p_a(B) + 1 \ \ \ &  \forall \ \ subcurve\ \ B\subseteq C \\
\  \deg \sG_{|B} \geq \ p_a(B) \ \ \ &  \forall \ subcurve\ \ B\subseteq C
\end{array}
$$

Then for every $n \geq 1$  the natural multiplication map $  (H^0(  C, \sH)) ^{\otimes n}
\rightarrow H^0( C,  \sH^{\otimes n})$ is surjective. 
\end{TEO} 

Moreover,  applying the same arguments used in  \cite[Proof of Thm. A, p. 327]{Fr} we have 

 \begin{PROP}
\label{teo:nef} 
 Let  $C$   be a   curve 
contained in a smooth algebraic surface
and
let  $\sH_1$, $\sH_2$ be two invertible sheaves such that 
 $\sH_1 \numeq \sF\otimes \sG_1$,    $\sH_2\numeq \sF \otimes \sG_2 $
  with 
 $$ \begin{array}{ll}
\  \deg \sF_{|B} \geq \ p_a(B) + 1 \ \ \ &  \forall \ \ subcurve\ \ B\subseteq C \\
\  \deg {\sG_1}_{|B} \geq \ p_a(B) \ \ \ & \forall \ subcurve\ \ B\subseteq C
\\ 
\  \deg {\sG_2}_{|B} \geq \ p_a(B) \ \ \ & \forall \ subcurve\ \ B\subseteq C
\\ 
\end{array}
$$ 
  (this holds e.g. if 
  $\sH_2 \numeq \sH_1\otimes \sA$,
  with $\sA$ a NEF invertible sheaf ). 
 
Then $  H^0(  C, \sH_1)  \otimes H^0(  C, \sH_2) \onto H^0(  C, \sH_1 \otimes  \sH_2) .$
\end{PROP} 

\section{ Even curves and even divisors}
	Let $C=\sum_{i=1}^{s} n_{i} \Gamma_{i}$ 
	be a 
	 curve 
	contained in a smooth algebraic surface.  
	
	 $C$ is said to be 
	{\em even} if 
	$$   \deg_{ \Gamma_{i}}  {\om_{C}}
	\mbox{ is even    for every  irreducible }  \ \Gamma_{i} \subset C $$
(see \S 1).  Similarly, if  $\sH$  is  an
	invertible sheaf on $C$, then $\sH$ is said to be 
	{\em even} if 
	$$  \deg_{ \Gamma_{i}}  {\sH}
\mbox{ is even  \ for every  irreducible } \ \Gamma_{i} \subset C $$
For even  invertible sheaves  of high degree the 
normal generation follows easily:
\begin{TEO}\label{teo:2}
 Let  $C=\sum_{i=1}^{s} n_{i} \Gamma_{i}$   be a  curve 
contained in a smooth algebraic surface
 and
let 
 $\sH $  be an even invertible sheaf on $C$ 
  such that 
  \begin{equation}\nonumber
	\  \deg_{B} \sH \geq \ 2 p_a(B) +2 \ \ \ \forall \ subcurve\ \ B\subseteq C
\nonumber
\end{equation}
Then for every $n \geq 1$  the natural multiplication map $  (H^0(  C, \sH)) ^{\otimes n}
\rightarrow H^0( C,  \sH^{\otimes n})$ is surjective. 
\end{TEO}

{\bf Proof.}	 First of all notice that $\sH$ is very ample by \cite[Thm. 1.1]{CFHR}. Moreover  
	since $\sH  $  is  even   there exists  an  invertible 
	sheaf $\sF$ such that ${\sF}^{\otimes 2} \numeq \sH$. By our numerical assumptions for every subcurve
	$B\subseteq C$
 it is 
	$  \deg_{B} \sF \geq \  p_a(B) +1 $ and $  \deg_{B}  (\sH\otimes \sF^{-1})\geq \  p_a(B) +1 $, whence 	
  we can conclude  by Thm.  \ref{teo:adjoint}.   \hfill 	  $\Box$

\section{The canonical ring of a  binary curve  }
In this section we deal with the particular case of binary curves, which is particularly interesting
from our point of view since it shows many of  the pathologies that may occur  studying 
curves with many components.
\begin{DEF} A curve $C$ is said to be a {\em binary curve}  if $C = \Gamma_{1}+ 
\Gamma_{2} $ 
with $\Gamma_{1}, \Gamma_2 $  irreducible and reduced rational curves such that 
$\Gamma_{1}
 \cdot \Gamma_{2}  = p_a(C) + 1 	$ 
(see  {\em \cite{cap}} )
 \end{DEF} 
If $C$ is not honestly hyperelliptic and $p_a(C) \geq 3$, then  
	   $\om_{C}$ is very ample on $C$
	by  \cite[Thm. 
	   3.6]{CFHR} and  
 $\varphi_{|\om_{C}|}$ embeds $C$ as the union of two 
rational normal curves  intersecting in a 0-dimensional scheme of 
length = $  p_a(C) + 1$. 

The theory of Koszul cohomology groups allow us to  give a deep analysis of the ideal ring of the embedded curves $\varphi_{|\om_{C}|}(C) \subset \proj ^{p_a(C) -1}$. 
The first proposition, which we will use in the proof of our main theorem,  is that 
  $R(C, \omega_C)$ is generated in degree 1. 
We point out that in \cite[Prop. 3]{CCM} there is an alternative proof of this result, but we prefer  to keep ours 
in order to show how the  theory of Kuszul cohomology groups  works.

\begin{PROP}\label{2components} 
	Let $C$ be a binary curve of genus $\geq 3$, 
 not honestly hyperelliptic. Then  $R(C, \omega_C)$ is generated in degree 1.
 \end{PROP}
{\bf Proof.} The proposition  follows if we 
 prove the vanishing of the Koszul groups 
$\sK_{0,q}(C, \om_C) $ $\forall \ q \geq 1 $.  

  For simplicity let 
$r = p_{a}(C) -1 = \deg \Oh_{\Gamma_i}( \omega_C)$  ($i=1,2$) and
 let us identify $\Gamma_{1} , \Gamma_{2} $ with their images in $\proj^r$. 
Notice that we have $H^{0}(C, \omega_{C}) \iso H^{0}(\Gamma_i, \omega_{C})$ (for $i=1,2$),  while it is
$H^0(\Gamma_1, \omega_C (-\Gamma_2))=0$ 
 and  $H^1(\Gamma_1, \omega_C^{\otimes k} (-\Gamma_2))=0$   for every $k \geq 2$. 
 Therefore  we can apply Prop. \ref{splitting} 
getting  the following exact sequence of Koszul groups 
$$ 
\cdots  \to  \sK_{1,q-1}(\Gamma_2, \om_C)    \to \sK_{0,q }(\Gamma_1,\Oh_{\Gamma_{1}}(-\Gamma_2),  \om_C)  \to \sK_{0,q}(C, \om_C ) \onto \sK_{0,q}(\Gamma_2, \om_C)   
$$
Now since $|\om_C|$ embeds $\Gamma_2$ as a rational normal curve in $\proj^r$ we have $ \sK_{0,q}(\Gamma_2, \om_C) =0$ for all $q \geq 1$. 
Moreover  by  \cite[(2.a.17)]{Gr}  we have 
$$ \sK_{0,q }(\Gamma_1,\Oh_{\Gamma_{1}}(-\Gamma_2),  \om_C)
\iso
\sK_{0,0 }(\Gamma_1, \om_C^{\otimes q} \otimes \Oh_{\Gamma_{1}}(-\Gamma_2),  \om_C) $$ 
But  on a rational curve we can easily determine the vanishing of the Koszul groups.
 Indeed,  by \cite[Corollary 3.a.6]{Gr},   if
$\Gamma_1 \iso \proj^1$, $\deg (\sF) = d , \deg (\sH) = a$ we have 
$$\sK_{p,0}( \Gamma_1, \sH, \sF)=0 \ \ \mbox{unless}   \ \ 0\leq a \leq 2d - 2 \ \mbox{and}  \ \ a-d+1 \leq p \leq a.  $$
Therefore  a simple  degree computation  yields   $ \sK_{0,q }(\Gamma_1,\Oh_{\Gamma_{1}}(-\Gamma_2),  \om_C)
=0 $ for $q \geq 3$, which implies
$\sK_{0,q}(C, \om_C) =0$ $\forall \ q \neq 2  $ by the above exact sequence of Koszul groups. 
If $q=2$,   it is 
$\sK_{0,2} (\Gamma_1, \Oh_{\Gamma_1}(-\Gamma_2 ), \omega_{C}) \iso \sK_{0,0 }(\Gamma_1, \om_C^{\otimes 2} \otimes \Oh_{\Gamma_{1}}(-\Gamma_2),  \om_C) \iso 
H^0( \Gamma_1 , \om_C^{\otimes 2}  \Oh_{\Gamma_{1}}(-\Gamma_2))\iso
H^0( \Gamma_1 , \om_C \otimes \om_{\Gamma_{1}})\iso
H^0( \proj^1 ,\Oh_{\proj^1} (r-2)).$
In particular  
  it is $\neq 0$.
  
   In this case, to show that  $\sK_{0,2}(C, \om_C) =0$ 
we are going to prove that the map $P: \sK_{1,1}(\Gamma_2, \om_C)   \to 
 \sK_{0,2 }(\Gamma_1,\Oh_{\Gamma_{1}}(-\Gamma_2),  \om_C)$ is surjective. 

Notice that 
it is  $ \sK_{1,1}(\Gamma_2, \om_C) \iso I_2(\Gamma_2, \proj^r )$, the space of quadrics 
in $\proj^r$ vanishing on $\Gamma_2$, and it is well known that  
$ I_2(\Gamma_2,  \proj^r )\iso \Ka^{  { r \choose 2 }  } $ (see e.g.,  \cite[Thm. 3.c.6]{Gr}).  
The idea of the proof  is   to show that  we can find $r-1$ points on $\Gamma_1$ 
  imposing  independent conditions on  
  the quadrics in 
$ I_2(\Gamma_2, \proj^{r})$. 
The proof will be  given  by induction on $r$.  	 
	
	 	 For $r=2$  it is 
		 $\sK_{1,1}(\Gamma_2, \om_C)   \onto  \sK_{0,2 }(\Gamma_1,\Oh_{\Gamma_{1}}(-\Gamma_2) )$
		since $C$ is embedded as a plane quartic.	
	For   $r \geq 3$   take    a point  
	 $Q \in \Gamma_1 \cap \Gamma_2$,  
	 	 choose coordinates $(x_0: ..:x_r) $  so that $Q=(0:...:0:1)$  and consider the projection from  $Q$
	 onto the hyperplane $\{x_r=0 \}$:
	$$ 
	\begin{array}{rcl}
	 \pi_Q: \proj^r  & \to & \proj^{r-1} \\
	 (x_0:...:x_r)& \mapsto & (x_0:...:x_{r-1})
	 \end{array} 
	 $$
	$\pi_Q (C) := \tilde{C} \subset \proj^{r-1}$ is again a binary curve
		canonically embedded in
	$ \proj^{r-1}$ by the complete system $| \omega_{\tilde{C}}|$.
	The projection map restricted to  $C$ and to its subcurves $\Gamma_1, \Gamma_2$
	  induces via pull-back  maps  between Koszul groups 
	of the same degree, which fit in the following commutative diagram  
	  \[ 
\xymatrix{ 
\sK_{1,1} (\tilde{C}, \omega_{\tilde{C}} )  \ar@{^{(}->}[d]^{\pi_Q^{\ast}}  \ar[r] 
 & \sK_{1,1} ({\Gamma_2} , \omega_{\tilde{C}} )  \ar@{^{(}->}[d]^{\pi_{2,Q}^{\ast}}
  \ar@{->}[r]^{ \hspace{ - 8mm}  P_{1}}  
 & \sK_{0,2} ({\Gamma_1}, \Oh_{\Gamma_1}(-{\Gamma_2} ), \omega_{\tilde{C}})  
  \ar@{^{(}->}[d]^{\pi_{1,Q}^{\ast}} 
 \\
\sK_{1,1} (C, \omega_{C} )  \ar@{->>}[d] \ar[r] 
 & \sK_{1,1} ( \Gamma_2 , \omega_C )  \ar@{->>}[d]\ar[r]^{ \hspace{ - 6mm} P} 
 & \sK_{0,2} (\Gamma_1, \Oh_{\Gamma_1}(-\Gamma_2 ), \omega_{C})   \ar@{->>}[d] 
 \\
\sQ (\omega_{C})   \ar@{^{(}->}[r] 
 & {\sQ}_2 (\omega_{C})    \ar[r]^{ \hspace{ - 6mm}  P_{2}}
 & \sR
}
\]
where 
$\sQ(\omega_{C})$  (resp. 
 ${\sQ}_2 (\omega_{C})$, $\sR$) denotes the cokernel of  $\pi_Q^{\ast}$ 
 (resp. the cokernel of ${\pi_{2,Q}}^{\ast}$, ${\pi_{1,Q}}^{\ast}$), and where   we 
have	  identified  the curves  $\Gamma_1$  and  $\Gamma_2$ with  their
 embeddings.

By our construction  ${\pi_{2,Q}^{\ast}} (\sK_{1,1} ({\Gamma_2} , \omega_{\tilde{C}} ))
\subset \sK_{1,1} ( \Gamma_2 , \omega_C ) $   is  the subspace spanned by the equations of the
 cones with
 center $Q$ over the  corresponding  quadric in $  I_2(\Gamma_2, \proj^{r-1}) 
  \iso \sK_{1,1} ({\Gamma_2} , \omega_{\tilde{C}}) $. Moreover   
  ${\sQ}_2 (\omega_{C})$ 
  is generated by  the following equations
 \begin{equation}\nonumber
 x_ix_r - x_{i+1}x_{r-1} \hspace{ 1 cm} \mbox{ for } \  i=0,..., r-2. 
 \end{equation} 
 (see  e.g    \cite[Ch. 6 , Prop. 6.1]{eis}).

To conclude  
  we  choose a 0-dimensional scheme $A$ of length $=r-1$ 
  so that $A'= A \cap \{x_r=0\}$ is a scheme 
  of length  $=r-2$ and  $\sR \iso \sI_{A'}/\sI_{A} \iso \Ka_{Q'} $ with  $Q'  \in \Gamma_1$
  in general position.  Since $\Gamma_1$ is rational it is immediately seen that 
   $$\sK_{0,2} (\Gamma_1, \Oh_{\Gamma_1}(-\Gamma_2 ), \omega_{C}) \iso 
H^0( \Gamma_1 , \om_C \otimes \om_{\Gamma_{1}}) \iso \Oh_A$$ and similarly 
$\sK_{0,2} ({\Gamma_1}, \Oh_{\Gamma_1}(-{\Gamma_2} ), \omega_{\tilde{C}})   \iso 
\Oh_{A'}.$
Now, by induction hypothesis   $P_1$  is  surjective, as well as 
   $P_2$  since $Q'  $
  is  general.  

 Therefore we can conclude that  $P: \sK_{1,1}(\Gamma_2, \om_C)   \to  \sK_{0,2 }(\Gamma_1,\Oh_{\Gamma_{1}}(-\Gamma_2),\omega_{C}) $ is onto, which is what we wanted to prove.
 \hfill$\Box$

	\section{Disconnecting components of numerically connected curves}

Taking an an irreducible component $\Gamma \subset C$  one  problem is that the restriction map
$H^0(C, \omega_C) \rightarrow H^0(\Gamma, {\omega_C}_{|\Gamma})$ is not surjective if 
$h^0(C-\Gamma, \Oh_{C-\Gamma}) =h^1(  C-\Gamma, \om_{C-\Gamma} )\geq 2$. 
Nevertheless,  if there exists  a curve $\Gamma$ with this property, it   plays a special role   in the proof
of  our main result.  

To be more explicit,  let us firstly consider the natural notion of 
disconnecting subcurve. 

\begin{DEF}
Let $C =\sum_{i=1}^{s} n_{i} \Gamma_{i}$  be a numerically connected curve. A subcurve $B \subset C$
 is said to be a {\em 
disconnecting subcurve}  if  
$h^0(C-B, \Oh_{C-B}) \geq 2$. 
\end{DEF}
If $B$ is a disconnecting curve then by the exact sequence
\[ H^0(C-B, \om_{C-B} )\to  H^0( C,  \om_{C}) \to 
 H^0( B,  \om_{C}) \to H^1(  C-B, \om_{C-B} ) \to  
 H^1( C,  \om_{C}) \] 
we deduce that  the restriction map 
$H^0(C, \omega_C) \rightarrow H^0(B, {\omega_C}_{|B})$ can not be surjective. In this case
 following the arguments
 pointed out in  \cite{Lau2} and \cite{Kon}  one can consider an ``intermediate" curve 
 $G$  such that  $B \subseteq G \subseteq C$ and
 $H^0(C, {\omega_C}) \onto H^0(G, {\omega_C}_{| G}) $.   
  We have the following useful Lemma. 

\begin{LEM}\label{disconnecting}
 Let $C =\sum_{i=1}^{s} n_{i} \Gamma_{i}$  be a $m$-connected curve  ($m \geq 1$)
 and $\Gamma \subset C$ be an irreducible and reduced disconnecting subcurve. 
 Let $G$ be a minimal subcurve of $C$ such that 
 $H^0(C, \omega_C) \onto H^0(G, {\omega_C}_{| G}) $  and  $\Gamma \subseteq G \subseteq C$.
 
  Setting  $E := C-G$, $G':= G-\Gamma$, 
 then 
 \begin{enumerate}
 \renewcommand\labelenumi{(\alph{enumi})}
 
  \item   $E$ is a maximal subcurve of  \  $C-\Gamma$ such that $h^1(E, \omega_E)=h^0(E, \Oh_E)=1$; 
 
  \item $\Gamma$ is of multiplicity $1$ in $G$, $\omega_G \otimes (\omega_C)^{-1}\iso \Oh_G(- E)$
  is NEF on  $G' $; 

  \item $\deg_{\Gamma}(E) = 
  \deg_{G'}(- E) + e $ with $e \geq m$;

  \item $h^1(E+\Gamma, \omega_{E+\Gamma})=1$, hence $H^0(C, \omega_C) \onto H^0(G',{\omega_C}_{|G'}) $;

   \item  
 $G$ is $m$-connected and in particular it is $h^1(G, \omega_G) =1$;

\end{enumerate}
\end{LEM}
{\bf Proof.} 
  By  hypotheses  $H^0(C, \omega_C) \not\onto H^0(\Gamma,{\omega_C}_{|\Gamma}) $ 
and $G$ is a  minimal subcurve such that $H^0(C, \omega_C) \onto H^0(G,{\omega_C}_{| G}) $. Therefore
$E = C-G $ is a maximal subcurve of $C-\Gamma$ 
such that $h^1(E, \omega_E)=h^1(C, \omega_C)=1$, proving ({\em a}). 

Moreover 
  by 
 \cite[Lemma  2.2.1]{Kon}
 either $\omega_G \otimes (\omega_C)^{-1}$ is NEF on $G$,  or $\Gamma$ is of multiplicity one in $G$ 
 and $\omega_G \otimes (\omega_C)^{-1}$ is NEF on $G-\Gamma =G'$. 
 
 Now  by adjunction it is $\omega_G \otimes (\omega_C)^{-1} \iso \Oh_G(- E)$, which has 
 negative degree on $G$   since $C$ is 
numerically connected by assumption. Therefore we can exclude the first case and by
\cite[Lemma  2.2.1]{Kon} 
we  conclude that  
$\Gamma$ is of multiplicity one in $G$, 
$\omega_G \otimes (\omega_C)^{-1} \iso \Oh_G(- E)$ is NEF on $G':=G-\Gamma $,  and  
$\deg_{\Gamma}(E) = 
  \deg_{G'}(- E) + e $ with $e \geq m$, proving $(b)$ and $(c)$.  
  
 To prove  $(d)$  consider the two curves $E$ and $\Gamma$. 
 It is 
  $h^1(E, \omega_E)=1$ by $(a)$ and  $h^1(\Gamma, \omega_{E+\Gamma})=0$
  because $\deg_{\Gamma}(\omega_{E+\Gamma} ) \geq 2p_a(\Gamma) -1$, whence we conclude
   considering the exact sequence
 \[ H^1( E, \omega_E) \rightarrow H^1(E+\Gamma, \omega_{E+\Gamma})
 \rightarrow  H^1(\Gamma, \omega_{E+\Gamma})=0\]
 
 $(e)$ follows 
since $\Oh_{G'}(-E)$ is NEF on $G'$. In fact  if $B \subset G$  without loss of generality 
we may assume $B \subset G'$,   and then we obtain  
$B\cdot (G-B)= B\cdot (C-B) -E\cdot B \geq B\cdot (C-B)\geq m$   since it is $G=C-E$, that is $G$ is m-connected.
$h^1(G, \omega_G)=1$ follows by \cite[Thm. 3.3]{CFHR}.

\hfill $\Box$

\hfill\break
With an abuse of notation we will call a subcurve $E\subset C$ as in Lemma \ref{disconnecting} a {
\em maximal connected subcurve} of $C-\Gamma$.

\hfill\break
The above Lemma allows us to consider the splitting $C= G+E$ since by connectedness  both the restriction maps
$H^0(C, \omega_C) \rightarrow H^0(G,{\omega_C}_{|G})$ and $H^0(C, \omega_C) \rightarrow H^0(E,{\omega_C}_{|E})$
are surjective. 

Concerning the subcurve $G$ we have the following theorem.

\begin{TEO}\label{disconnecting-finale}
 Let $C =\sum_{i=1}^{s} n_{i} \Gamma_{i}$  be an even  $4$-connected curve and assume there exists  
 an irreducible and reduced disconnecting subcurve $\Gamma \subset C$. 	  
 
  Let $G$ be a minimal subcurve of $C$ such that 
 $H^0(C, \omega_C) \onto H^0(G, \omega_C) $  and  $\Gamma \subseteq G \subseteq C$. 
 
 Then on $G$ the
	  multiplication map 
	 $H^0(G, \omega_C) \otimes H^0(G,\omega_G) \to H^0(G, \omega_C \otimes \omega_G) $ is surjective.
\end{TEO}

 To simplify the notation, for every subcurve 
  $B\subset C$ by  $H^0(B, \omega_C)$ we will denote 
 the space of sections of  $ {\omega_C}_{|B}$.
 
 \hfill\break If there exists a disconnecting component $\Gamma$ and a   decomposition $C=G'+\Gamma +E$ as in  Lemma \ref{disconnecting}  such that 
$h^1(G', \omega_{G'})\geq 2 $ then we need an auxiliary Lemma. 

\begin{LEM}\label{disconnecting3a}
 Let $C =\sum_{i=1}^{s} n_{i} \Gamma_{i}$  be an even  $4$-connected curve and assume there exists  
 an irreducible and reduced disconnecting subcurve $\Gamma \subset C$. 

If there exists a decomposition  $C=G'+\Gamma +E$ as in  Lemma \ref{disconnecting}  such that 
$h^1(G', \omega_{G'})\geq 2 $ 
  then  there exist a decomposition $C=E+ \Gamma +G_1+G_2$ s.t. 

\begin{enumerate}
 \renewcommand\labelenumi{(\alph{enumi})}
 
  \item  $G_2 + \Gamma$ is 4-connected 
  
 \item  $h^1(G_1, \omega_{G_1}) =1$ 
   
\item  $\Oh_{G_2}(- G_1 )$ is NEF on $G_{2} $

 \item $H^0(G, \omega_G) \onto H^0(G_2+\Gamma, \omega_G) $.

\end{enumerate}

\end{LEM}
{\bf Proof.}  Let 
 $C= E + G$ and  $G= \Gamma + G'$ be as in Lemma \ref{disconnecting}. By ({\em e}) of Lemma \ref{disconnecting} $G$ is 4-connected and by our hypothesis $h^1(G', \omega_{G'}) \geq 2$, i.e.,  the irreducible curve $\Gamma$
is a disconnecting component for $G$ too. Therefore by Lemma \ref{disconnecting} applied to $G$, there exists a maximal connected subcurve $G_1 \subset G'$  
 and a decomposition 
 $G= \Gamma+G_1 +G_2$ such that ({\em a}), ({\em b}), ({\em c}), ({\em d}) hold. 
\hfill $\Box$ 

\hfill\break
{\bf Proof of Thm. \ref{disconnecting-finale}.}   The proof of theorem \ref{disconnecting-finale} will be treated considering separately the case $h^1(G', \omega_{G'})=1$ and  $h^1(G', \omega_{G'})\geq 2.$

\hfill\break{\bf Case 1:} 
{\em There exists a disconnecting component $\Gamma$ and a decomposition $C=G'+\Gamma +E$ as in  Lemma \ref{disconnecting} such that 
 $h^1(G', \omega_{G'})=1$.}

 Let $G=G' + \Gamma$. On 
$\Gamma$ both the invertible sheaves 
$\om_{\Gamma}(E) $ and $\omega_G$ 
have degree $\geq 2p_a(\Gamma)+2$. 
In particular  we have the following 
 exact sequence 
\[
0\ra H^0(\Gamma, \om_{\Gamma}(E) ) \ra 
 H^0( G, \om_{C}) \ra
 H^0(G', \om_{C} ) \ra 0  
\] 
Twisting with $H^0(  \om_{G})= H^0( G, \om_{G})$ 
we get  the following commutative diagram:
\[
\begin{array}{ccccc}
 H^0(\Gamma, \om_{\Gamma}(E) ) \otimes H^0(  \om_{G}) &\into & 
 H^0( G, \om_{C}) \otimes H^0(  \om_{G}) &
\onto& H^0(G', \om_{C} ) \otimes H^0( \om_{G})   \\
  r_{1}  \downarrow \ \ \ & & r_{2} \downarrow \ \ \ &  & r_{3}\downarrow  
  \ \ \ \\ 
 H^0(\Gamma ,\om_{\Gamma} (E) \otimes \om_{G}) &\into & H^0( G, \omega_C \otimes \omega_G) &
\onto & H^0(G', \omega_C \otimes \omega_G)   \\
\end{array}
\] 
Now, since  by our hypothesis $h^1(G', \omega_{G'})=1$ then   $H^0(G,\omega_G)\onto H^0( \Gamma,  \om_{G})$  and 
 by   \cite[Thm.6] {Mu}  we have the surjection  
$H^0(\Gamma, \om_{\Gamma}(E) ) \otimes H^0( \Gamma,  \om_{G})  \onto 
 H^0(\Gamma ,\om_{\Gamma} (E) \otimes \om_{G}) $.

  The proposition follows since  also  $r_{3}$ is surjective by Prop. \ref{teo:nef}.  Indeed, 
  it is ${\omega_G}_{| G'} \iso {\omega_{C}}_{| G'} (-E)$  with $\Oh_{G'}(-E)$ NEF and 
  ${\omega_{C}}_{| G'}$ is an even  invertible sheaf whose degree on every subcurve $B\subseteq  G'$ 
  satisfies  $\deg_{B} (\omega_C) \geq 2p_a(B)+2$.

\hfill\break{\bf Case 2:} 
{\em There exists a disconnecting component $\Gamma$ and a decomposition $C=G'+\Gamma +E$ as in  Lemma \ref{disconnecting} such that 
 $h^1(G', \omega_{G'})\geq 2$.}
 
 Let $C=E+ \Gamma +G_1+G_2$ and $G=\Gamma +G_1+G_2$,  be a decomposition as in Lemma \ref{disconnecting3a}. 
We proceed as in Case 1,  considering the curve $ G_2+\Gamma$ instead of the irreducible $\Gamma$.

 First of all  let us prove that $H^1(G_2 + \Gamma, \om_{G_2+\Gamma}(E) )=0$. 
 
It is
  $ ( \om_{G_2+\Gamma}(E) )_{|G_2} \iso 
 (\omega_C (-G_1))_{|G_2}$ and  in particular 
for every subcurve   $B\subseteq G_2$ it is 
$\deg_{ B} (\om_{G_2+\Gamma}(E) )\geq 2p_a( B) +2$
 since $\Oh_{G_2}(-G_1)$ is NEF. 

 If   $B\not\subseteq G_2$, we can write  $B=B'+\Gamma$,  with $B'\subseteq G_2$,  obtaining   
 $$\deg_{ B} (\om_{G_2+\Gamma}(E))=
 \deg_{ B} (\om_{G_2+\Gamma})+ E\cdot  B \geq 2p_a( B) +2$$
since 
 $E\cdot  B=E\cdot (B' + \Gamma)  \geq E\cdot (G' + \Gamma) \geq 4 $ and $ \deg_{ B} (\om_{G_2+\Gamma}) \geq 2p_a( B) -2$
 by connectedness.
Therefore 
 $H^1(G_2 + \Gamma, \om_{G_2+\Gamma}(E) )=0$ by \cite[Lemma 2.1]{CF}
 and  we have the following exact sequence 
\[
0\ra H^0(G_2 + \Gamma, \om_{G_2+\Gamma}(E) ) \ra 
 H^0( G, \om_{C}) \ra
 H^0(G_1, \om_{C} ) \ra 0  
\] 
Twisting with $H^0(  \om_{G})= H^0( G, \om_{G})$  we can argue as in 
Case 1. Indeed, consider the 
	commutative diagram:
\[
\begin{array}{ccccc}
 H^0(G_2 +  \Gamma, \om_{G_2 + \Gamma}(E) ) \otimes H^0(  \om_{G}) &\into & 
 H^0( G, \om_{C}) \otimes H^0(  \om_{G}) &
\onto& H^0(G_1, \om_{C} ) \otimes H^0( \om_{G})   \\
  r_{1}  \downarrow \ \ \ & & r_{2} \downarrow \ \ \ &  & r_{3}\downarrow  
  \ \ \ \\ 
 H^0(G_2 +  \Gamma ,\om_{G_2 + \Gamma} (E) \otimes \om_{G}) &\into & H^0( G, \omega_C \otimes \omega_G) &
\onto & H^0(G_1, \omega_C \otimes \omega_G)   \\
\end{array}
\] 
The map $r_3$ is onto by Prop. \ref{teo:nef} since $\om_G \iso \om_C(-E)$,
$\Oh_{G_1}(-E)$ is NEF 
 and by 
Lemma \ref{disconnecting3a}  we have the surjection $H^0( \om_{G})  \onto H^0(G_1, \omega_G)$.

 The Theorem  follows if we show that  $r_1$ 
 is surjective too.  
 Notice that we can write the multiplication map $r_1$ as follows
 $$r_1:  H^0(G_2 +  \Gamma, \om_{G_2 + \Gamma}(E) ) \otimes  H^0(G_2 +  \Gamma, \om_{G_2 + \Gamma}(G_1) ) 
 \to  H^0(G_2 +  \Gamma ,\om_{G_2 + \Gamma}^{\otimes 2}  (E+G_1)) $$
that is, $r_1$ is symmetric in $E$ and $G_1$.

Assume firstly that $ (G_1-E)\cdot \Gamma  \geq 0$. We proceed  considering
 a general effective Cartier divisor  $\Upsilon$ on $G_2+\Gamma$ such that 
 \begin{equation}\nonumber  
	 \left\{ \begin{array}{rcll}
 ( \Oh_{G_2}(  \Upsilon) )^{\otimes 2} & \numeq&  {\omega_C (- 2E )}_{|G_2}  \\ 
  \deg( \Oh_{\Gamma}(  \Upsilon) ) &= & \frac{1}{2} \deg({\omega_C}_{|\Gamma}) -  E\cdot \Gamma - \delta   
 \end{array} \right.
		\end{equation} 
		with $\delta = \lceil  \frac{-G_1 \cdot G_2}{2} \rceil$.  
We remark that  by connectedness of $C$ and nefness of $\Oh_{G_2}(-2E)$
 it is $\deg (\Oh_{\Gamma_i}(\Upsilon)) \geq 1$ for every $\Gamma_i\subseteq G_2$ 
and  by our numerical conditions it is
$$ \deg( \Oh_{\Gamma}(  \Upsilon) ) 
		=  		 p_a(\Gamma) -1  +\frac{1}{2} (G_1 -E)\cdot \Gamma  + \frac{1}{2} G_2 \cdot  \Gamma -\delta \geq 1$$
since 
$ (G_1-E)\cdot \Gamma  \geq 0$
		by our assumptions and 
$G_2 \cdot  \Gamma   - 2\delta \geq (\Gamma+G_1 )\cdot G_2 \geq (\Gamma+G_1+ E )\cdot G_2  \geq 4$ by 4-connectedness of $C$.
		
Now let  
$\sF := \omega_G(-\Upsilon)$. $\sF$  
is a general invertible subsheaf of ${\omega_G}_{| G_2 +\Gamma}$ s.t.
 \begin{equation}\nonumber  
	 \left\{ \begin{array}{rcll}
		\deg_{B}  {\sF}   & = &\frac{1}{2}\deg_{B}{\omega_C} \ \ \forall B \subseteq G_2 \\ 
\deg_{\Gamma}  {\sF}  &= &  \frac{1}{2} \deg _{\Gamma}{\omega_C} + \delta  
	\end{array} \right.
		\end{equation} 
 By 4-connectedness  on every subcurve $B\subseteq  G_2+\Gamma$ it is  $\deg_{B}  {\sF} \geq p_a(B)+1$.
By Theorem  \ref{teo:general}  we conclude that  
 $|\sF|$ is base point free and  $h^1(G_2 +\Gamma, \sF)=0$. 
 
Therefore we have  the following exact sequence
$$
0
\to H^0(G_2 +  \Gamma, \sF)    \to  H^0(G_2 +  \Gamma, \omega_{G})  \to  H^0( \Oh_{  \Upsilon} )\to 0 
$$
and we have the surjection  $H^0( \Oh_{  \Upsilon}) \otimes H^0(G, \om_{G_2 + \Gamma}(E) ) \onto H^0(G, \Oh_{  \Upsilon} \otimes \om_{G_2 + \Gamma}(E) ) $
since $\Oh_{  \Upsilon}$ is a skyscraper sheaf and $| \om_{G_2 + \Gamma}(E) |$ is base point free by \cite[Thm.3.3]{CFHR}. 
Whence the map $r_1$ 
 is surjective if we  prove that  
\begin{equation}\label{onto}
\nonumber H^0(G_2 +  \Gamma, \om_{G_2 + \Gamma}(E) ) \otimes H^0(G_2 +  \Gamma,  \sF) \onto 
H^0(G_2 +  \Gamma, \om_{G_2 + \Gamma}(E) \otimes \sF)\end{equation}
To this aim we are going to apply  ({\em  iii}) of Prop. \ref{prop:van}.
First notice that
 $$ H^0(G_2 +  \Gamma, \sF) \into 
H^0( \Gamma, \sF)$$ 
Indeed, by adjunction and Serre duality the kernel of this map is isomorphic to 
$H^0(G_2,  \sF  - \Gamma)  \iso H^1(G_2, \omega_C(-E-G_1) \otimes {\sF}^{-1})$,  which   vanishes  by Thm. \ref{teo:general} 
since  it
is the first cohomology group of   a general  invertible sheaf whose degree  on every component  $B\subseteq G_2$ satisfies
 $$ \deg_{B}(\omega_C(-E-G_1) \otimes {\sF}^{-1})=\frac{1}{2}( \deg_{ B} \omega_C  ) + (-E-G_1)\cdot B \geq  \frac{1}{2}( \deg_{ B}  (\omega_C)) \geq p_a(B)$$
because  $C$ is 4-connected and $\Oh_{G_2}(-E-G_1)$ is NEF. 

Moreover we have also the embedding 
$$ H^0(G_2 +  \Gamma,  \om_{G_2 + \Gamma} \otimes [\om_{G_2 + \Gamma}(E)]^{-1} \otimes {\sF} ) \into
 H^0(\Gamma,  \om_{G_2 + \Gamma} \otimes [\om_{G_2 + \Gamma}(E)]^{-1} \otimes {\sF} )
$$
since  $H^0(G_2,  \sF - E - \Gamma) \iso H^1(G_2, \omega_C(-G_1) \otimes {\sF}^{-1})=0$
because 
 $$ \deg_{B}(\omega_C(-G_1) \otimes {\sF}^{-1})=\frac{1}{2}( \deg_{ B} \omega_C  ) + (-G_1\cdot B) \geq  \frac{1}{2}( \deg_{ B}  (\omega_C)) \geq p_a(B)$$
by  4-connectedness of $C$ and nefness of  $\Oh_{G_2}(-G_1)$.
 
In order to conclude we are left to compute $ h^0(G_2 +  \Gamma,  \om_{G_2 + \Gamma} \otimes [\om_{G_2 + \Gamma}(E)]^{-1} \otimes {\sF} )= h^0(G_2 +  \Gamma, \sF (-E) ).  $ 

 $ {\sF} (-E)$ 
 is a general invertible sheaf s.t.  
   \begin{equation}\nonumber  
	 \left\{ \begin{array}{rcll} 
	\deg_{B} ({\sF} (-E))   & = &\frac{1}{2}\deg_{B}\omega_C  -E \cdot B   \ \ \ \  \ \ \ 
\forall B \subseteq G_2   	  \\	
		\deg_{\Gamma} ( \sF (-E))  &= &  
		 \frac{1}{2}\deg_{\Gamma}(\omega_{G_2+\Gamma}) + \frac{1}{2}( G_1-E)\cdot\Gamma +\delta  
		 \\
	\end{array} \right.
		\end{equation} 
Therefore 
we obtain immediately that $\deg_{B} \sF (-E)  \geq p_a(B)$  on every $B\subseteq G_2$,  
whereas if $B=B'+\Gamma$ with $B'\subseteq G_2$ 
it is 
$$ \deg_{B} ({\sF} (-E))  =  \frac{1}{2}\deg_{B} (\omega_{G_2+\Gamma})+ \frac{1}{2}(G_1 -E)\cdot B' + 
 \frac{1}{2}( G_1-E)\cdot\Gamma +\delta \geq p_a(B)$$
  since 
by our assumptions $G_2+\Gamma$ is numerically connected,   $ (G_1-E)\cdot \Gamma  \geq 0$, $ E\cdot B'\leq 0$   and 
$\delta \geq \frac{1}{2}(-G_1\cdot G_2)\geq  \frac{1}{2}(-G_1\cdot B') $.

 In particular  
 by Theorem 
 \ref{teo:general} we get  
$H^1(G_2 +  \Gamma, \sF (-E))=0$ and 
 by Riemann-Roch theorem  we have $h^0(G_2+\Gamma, \sF(-E))= h^0(G_2+\Gamma, \sF) - E\cdot (\Gamma+G_2)$. 
Finally, since $\Oh_{G_{1}}(-E)$ is NEF we have  $E\cdot (\Gamma+G_2) 
 \geq  E\cdot (\Gamma+G_2+G_1) \geq  4$
 because $C$ is 4-connected, that is,  $h^0(G_2+\Gamma, \sF(-E))\leq  h^0(G_2+\Gamma, \sF) -4$.
Whence all the hypotheses of ({\em iii}) of Prop. \ref{prop:van} are satisfied and we can conclude.

If  $ (E - G_1)\cdot \Gamma  < 0$  we exchange the role of $\Oh_{G_2+\Gamma}(E)$ with the one of 
$\Oh_{G_2+\Gamma}(G_1)$ and we reply the proof ``verbatim", since our numerical conditions
are symmetric in $E$ and $G_1$.

\hfill $\Box$

	\section{The canonical ring of an even 4-connected curve }
In this section we are going to   show  Theorem \ref{teo:3}. 

\hfill\break
{\bf Proof of  Theorem \ref{teo:3}.} 
For all $k\in \Na$ 
we have to show the surjectivity of the maps
\[\rho_k: (H^{0}(C,\omega_C))^{\otimes k} \longrightarrow H^{0}(C,\omega_C^{\otimes k}) \]
 For $k=0,1$ it is obvious.  For $k\geq 3$ it follows by an induction argument applying  Prop. \ref{prop:van}
 to the sheaves
$\omega_C^{\otimes (k-1)}$ and $\omega_C$.

For $k=2$ the proof is  based on the above results. 
If $C$ is irreducible and reduced  the result is (almost) classical. 
For the general case 
we separate the proof in three different parts, 
depending on the existence of suitable irreducible components. 

\vspace{-2 mm}

\hfill\break
{\bf   Case A: } {\em There exists a  not disconnecting irreducible  curve $\Gamma$  of arithmetic genus $p_a(\Gamma)\geq 1$.}

\hfill\break
In this case, writing $C= \Gamma + E$, we have  the surjections 
$H^0(C, \omega_C) \onto H^0(\Gamma,{\omega_C}_{|\Gamma})$ and 
$H^0(C, \omega_C) \onto H^0(E,{\omega_C}_{|E})$, whence   we can conclude by 
 the following 
	commutative diagram:
\[
\begin{array}{ccccc}
 H^0(\Gamma, \om_{\Gamma} ) \otimes H^0(  \om_{C}) &\into & 
 H^0(  \om_{C}) \otimes H^0(  \om_{C}) &
\onto& H^0(E, \om_{C} ) \otimes H^0( \om_{C})   \\
  r_{1}  \downarrow \ \ \ & & \rho_{2} \downarrow \ \ \ &  & r_{3}\downarrow  
  \ \ \ \\ 
 H^0(\Gamma ,\om_{\Gamma}\otimes \om_{C}) &\into & H^0( C, \om_{C}^{\otimes 2}) &
\onto & H^0(E, \om_{C}^{\otimes 2})   \\
\end{array}
\] 
(where    $H^0(  \om_{C})= H^0( C, \om_{C})$).  Indeed, 
since $C$ is  $4$-connected and $\om_{C}$ is  an even divisor  we get 
the surjection of the map $r_{3}$  by  
 Theorem \ref{teo:2},  while  Proposition \ref{cor:irreducible}  
 ensure 
 the surjectivity  of the map  $r_{1}$, forcing $\rho_2$ to be surjective too 
 (cf.  also \cite[Thm. 6]{Mu}). 

\hfill\break
{\bf  Case B:}  {\em 
There exists a disconnecting irreducible component $\Gamma$. } 

\hfill\break
Let us consider the decomposition $C= E + G $ introduced in Lemma \ref{disconnecting}. Then we have
the exact sequence
$$ 
0\to H^0(G ,\om_{G}) \to  H^0( C, \om_{C}) 
\to  H^0(E, \om_{C}) \to 0    
$$ 
and furthermore by Lemma \ref{disconnecting} {\em (i)}  also the map  $H^0( C, \om_{C}) 
\to  H^0(G, \om_{C})$ is onto.    Replacing  $\Gamma$ with $G$ we can build  a commutative diagram analogous 
to the one shown in case A.  
Keeping the notation $r_1, r_3$ for the analogous  maps,   
   by 
 Theorem \ref{teo:2}  and Theorem \ref{disconnecting-finale} the maps $r_3$ and $r_1$  are surjectives,
whence also $\rho_2$ is onto.

\hfill\break
{\bf  Case C: } {\em Every irreducible component $ \Gamma_{i} $ of $C$ has  arithmetic genus 
$p_a(\Gamma_i)=0 $ and it 
 is not disconnecting. } 

\hfill\break
First of all notice that by 
connectedness  for every irreducible $\Gamma_h$ there exists at least one
$\Gamma_k$ such that $\Gamma_h \cdot \Gamma_k \geq 1$, and that $h^0(B,\Oh_B)=1$ for a curve $B=\sum a_i\Gamma_i \subset C$ implies $\Gamma_i \cdot (B-\Gamma_i) \geq 1$ for every $\Gamma_i\subset B$.

 We will consider separately the different situations that may happen.

 \hfill\break
{\bf {\em   C.1.} }
{\em There exist  two components $\Gamma_h, \Gamma_k$  (possibly $h=k$ if $\mult_{C}\Gamma_h \geq 2$) such that
 $\Gamma_h \cdot \Gamma_k \geq 2$,  and  
 $\Gamma= \Gamma_h +\Gamma_k \subset C$  is not disconnecting.}
 
  Notice that by \S 4 we may assume $C-\Gamma \neq \emptyset $.
 In this situation, setting $E=C-\Gamma$ 
by  {\em (ii)} of Proposition \ref{cor:irreducible} we can proceed exactly as in {\em Case A}. 
   
   \hfill\break
{\bf {\em   C.2.} }  {\em  There exist  two components $\Gamma_h, \Gamma_k$ (possibly $h=k$ if $\mult_{C}\Gamma_h \geq 2$)  such that 
 $\Gamma:= \Gamma_h +\Gamma_k \subset C$  is disconnecting and  $\Gamma_h\cdot \Gamma_k \geq 0$. }
 
In this case take $E$ a maximal subcurve of 
$C-\Gamma_k- \Gamma_h$ such that $h^0(E, \Oh_E)=1$ and let $G=C-E$.   Then  we obtain a decomposition 
$C=E+\Gamma +G' $ with $\Oh_{G'}(-E)$  NEF on $G'$.

Firstly let us point out some useful  remarks about this decompositions. 

We have   $h^0(E+\Gamma_k+G', \Oh_{E+\Gamma_k+G'})=1$ 
since 
$\Gamma_k$  is  not  disconnecting in  $C$ and  it is immediately seen that also 
$ h^0(\Gamma_k+G', \Oh_{\Gamma_h+G'})=1$ because  $\Oh_{G'}(-E)$  is NEF on $G'$. 
But  
 $p_a(\Gamma_k)=0$,   whence by the remark given at the beginning of Case C it is  $\Gamma_k \cdot G' \geq 1$.  In particular 
 $H^0(G,\omega_G) \onto H^0(\Gamma_h,\omega_G) $. 
 Similarly  we obtain  $\Gamma_h \cdot G' \geq 1$, $H^0(G,\omega_G) \onto H^0(\Gamma_k,\omega_G) $,  and considering
 $\Oh_{\Gamma}(E)$ we have 
$E\cdot \Gamma_h \geq 1$  and $E\cdot \Gamma_k \geq 1$. 
 Furthermore, since $E \cdot \Gamma \geq 4$, may assume $E\cdot \Gamma_h\geq 2$ . \\

We will consider firstly the subcase where  $\Gamma_h\cdot \Gamma_k \geq 1$ and secondly the case where the product is null.

  \hfill\break
{\bf {\em   C.2.1.} }     If   $\Gamma_h\cdot \Gamma_k \geq 1$ and $\Gamma=\Gamma_h +\Gamma_k \subset C$  is disconnecting, 
arguing as in  {\em  Case B},  the theorem  follows if  we have the surjection of the multiplication map
 $r_1:H^0(G,{\omega_C}) \otimes H^0(G,\omega_G) \to H^0(G, \omega_C \otimes \omega_G) $.
Considering  the diagram
\[
\begin{array}{ccccc}
 H^0(\Gamma, \om_{\Gamma}(E) ) \otimes H^0(  \om_{G}) &\into & 
 H^0( G, {\om_{C}}) \otimes H^0(  \om_{G}) &
\onto& H^0(G',  {\om_{C}}) \otimes H^0( \om_{G})   \\
  s_{1}  \downarrow \ \ \ & & r_1  \downarrow \ \ \ &  & t_{1}\downarrow  
  \ \ \ \\ 
 H^0(\Gamma ,\om_{\Gamma} (E) \otimes \om_{G}) &\into & H^0( G, \omega_C \otimes \omega_G) &
\onto & H^0(G', \omega_C \otimes \omega_G)   \\
\end{array}
\] 
 it is sufficient to show that  $s_1$ is onto since $t_{1}$ is surjective by Prop. \ref{teo:nef}. 
 To this aim we take  the splitting 
 $$ 0\to H^0(\Gamma_h, \omega_{\Gamma_h} (E)) \to
 \ H^0(\Gamma, \omega_{\Gamma} (E)) \to \ H^0(\Gamma_k, \omega_{\Gamma_k} (\Gamma_h +E)) \to 0 
 $$
 Twisting with $H^0(G, \omega_G) = H^0(\omega_G)$  (notice that  $H^0(G,\omega_G) \onto H^0(\Gamma_h,\omega_G) $ 
and similarly for $\Gamma_k$  by the above remark),
  we can conclude since we have the surjections 
 
 $$ H^0(\Gamma_h, \omega_{\Gamma_h} (E)) \otimes H^0(\Gamma_h, \omega_G) \onto 
 H^0(\Gamma_h, \omega_{\Gamma_h} (E)\otimes \omega_G)$$ 
 $$H^0(\Gamma_k, \omega_{\Gamma_k} (\Gamma_h +E)) \otimes H^0(\Gamma_k, \omega_G) \onto 
H^0(\Gamma_k, \omega_{\Gamma_k} (\Gamma_h +E)\otimes \omega_G)$$ 
 because     $\Gamma_h \iso \Gamma_k \iso \proj^1$,   and  all the sheaves have  positive degree on both the curves
 (see \cite[Corollary 3.a.6]{Gr} for details).

 \hfill\break
{\bf {\em   C.2.2.} } Assume now   $\Gamma_h\cdot \Gamma_k =0$ and $\Gamma=\Gamma_h +\Gamma_k \subset C$  to be disconnecting.

  If  $E\cdot \Gamma_h\geq 2$, $E\cdot \Gamma_k\geq 2$ and $G'\cdot \Gamma_h\geq 2$, $G'\cdot \Gamma_k\geq 2$  then 
  we consider  the exact sequence 
  $$0 \to \om_{\Gamma}(E)  \to {\omega_C}_{| (G'+\Gamma)}  \to {\omega_C}_{| G'}  \to 0$$ 
  and we operate as in {\em   C.2.1.} 
   
   Otherwise, without loss of generality, we may assume $E\cdot \Gamma_h=1$ or $G' \cdot \Gamma_h=1$. \\ 
   
 If $E\cdot \Gamma_h=1$  
   then by 4-connectedness of $C$ it is $E\cdot \Gamma_k\geq 3$,
    $G' \cdot   \Gamma_h\geq 3$,  $G' \cdot   \Gamma_k\geq 3$.
    
    Let $G= G'+  \Gamma_h +  \Gamma_k $ and consider the splitting $C=E+G$: 
 as in the previous case it is enough
    to prove the surjection of 
 $r_1:H^0(G,{\omega_C}) \otimes H^0(G,\omega_G) \to H^0(G, \omega_C \otimes \omega_G) $. 
 
 To this aim we take the following exact sequence   
    \[
0 \to %
 H^0(\Gamma_k, \om_{\Gamma_k}(E) )  \to 
 H^0( G, {\om_{C}})\to 
H^0(G' +\Gamma_h,  {\om_{C}})\]
  By our numerical conditions $|\omega_G|$ is base point free and by connectedness $H^0(\omega_G)\onto H^0(\Gamma_k,\omega_G).$ 
    
By  \cite[(2.a.17), (3.a.6)]{Gr} (or simply since we have sheaves of positive degree on a rational curve)  we have the surjection  $H^0(\Gamma_k,\omega_G) \otimes H^0(\Gamma_k,\omega_{\Gamma_k}(E))
 \onto H^0(\Gamma_k,\omega_G \otimes \omega_{\Gamma_k}(E))$.

On the contrary it is $\deg\Oh_{G'+\Gamma_h}(-E)\geq -1$.   
    Therefore  we can consider a subsheaf  $\sF \subset {\omega_C}_{|G'+\Gamma_h}$  such that $(\sF_{|G'+\Gamma_h})^{\otimes 2} \numeq {\omega_C}_{|G'+\Gamma_h}$.
    Then  for every
    $B\subset G'+\Gamma_h$ 
$\sF_{|B}$ has degree at least $p_a(B)+1$   whilst 
  $\omega_G \otimes \sF^{-1}$ is an invertible sheaf  of degree at least $p_a(B)$.    Whence  by Prop. 
    \ref{teo:nef}  $H^0(G' +\Gamma_h , \omega_{G}) \otimes H^0(G'+\Gamma_h , {\omega_{C}}) \onto
     H^0(G'+\Gamma_h, \omega_{G}
  \otimes{\omega_{C}})$ and then $r_1$ is onto. \\

     If $G'\cdot \Gamma_h=1$ then by 4-connectedness of $C$ it is $E\cdot \Gamma_h\geq 3$,
    $E \cdot   \Gamma_k \geq 3$,  $G' \cdot   \Gamma_k\geq 3$.

  In this case we   
     write $\tilde{E}:= E+\Gamma_h$,  $\tilde{G}:= G' +\Gamma_k$.  We have a decomposition $C= \tilde{E} + \tilde{G}$  where 
  $\tilde{E}$ is connected, $\tilde{G} $ is  3-connected. Moreover 
  by adjunction we have the isomorphism $\omega_{\tilde{G}}=  {\omega_{C} (-\tilde{E}) }_{|\tilde{G}}$, where
$ \deg \Oh_{\tilde{G}}(-\tilde{E})\geq -1$. 

Arguing as above the theorem follows  if we 
  prove the surjection of 
 $\tilde{r}_1:H^0(\tilde{G},{\omega_C}) \otimes H^0(\tilde{G},\omega_{\tilde{G}}) \to H^0(G, \omega_C \otimes \omega_{\tilde{G}}) $. 

To this aim let us show firstly that $h^0(G', \Oh_{G'})=1$. Indeed, 
since  $\Gamma_k$ is not disconnecting for $C$, whilst it is 
disconnecting for $C-\Gamma_h$  it is $h^0(G'+\Gamma_h, \Oh_{G'+\Gamma_h})=1$. Therefore 
since $\deg \Oh_{\Gamma_h} (-G') = -1$  and $\Gamma_h \iso \proj^1$ by 
 the exact sequence
\[0= %
 H^0(\Gamma_h, \Oh_{\Gamma_h} (-G') ) \to 
 H^0(G'+\Gamma_h, \Oh_{G'+\Gamma_h})\to 
H^0(G', \Oh_{G'}) \to  H^1(\Gamma_h, \Oh_{\Gamma_h} (-G'))=0\]
we obtain $h^0(G', \Oh_{G'})=1$.
Going back to $\tilde{G}:= G' +\Gamma_k$  let us take  the   exact sequence   
    \[
0 \to %
 H^0(\Gamma_k, \om_{\Gamma_k}(E) )  \to 
 H^0( \tilde{G}, {\om_{C}})\to 
H^0(G' ,  {\om_{C}})\to 0. \]
We have $H^0(\Gamma_k, \om_{\Gamma_k}(E) ) \otimes H^0(\Gamma_k, \omega_{\tilde{G}}) \onto 
H^0(\Gamma_k, \om_{\Gamma_k}(E)  \otimes \omega_{\tilde{G}})$  since $\Gamma_k \iso \proj^1$.

Finally,  taking a  subsheaf  $\sF \subset {\omega_C}_{|{\tilde{G}}}$  such that $(\sF_{{\tilde{G}}})^{\otimes 2} \numeq {\omega_C}_{|{\tilde{G}}}$, we can consider the two sheaves 
   $\sG_1:=\omega_C\otimes \sF^{-1} \numeq \sF$,  
  $\sG_2:=\omega_{\tilde{G}} \otimes \sF^{-1}$. 
  
  $\sF, \sG_1, \sG_2$  satisfy the assumptions of   Prop. 
    \ref{teo:nef},  whence         $\tilde{r}_1$ is surjective  and we can conclude.

 \hfill\break
 {\bf  {\em   C.3.}  } 
{\em    There exists  two distinct irreducible components 
 $\Gamma_h, $ $\Gamma_k$ such that 
 $\Gamma_h\cdot \Gamma_k=0$,
  and  
 $\Gamma= \Gamma_h + \Gamma_k$ is 
not disconnecting. }

In this case 
   the situation is slightly different.
   
    By \S 2.3 
  $\rho_2: (H^{0}(C,\omega_C))^{\otimes 2} \onto H^{0}(C,\omega_C^{\otimes 2})$ iff 
  $\sK_{0,1}(C,\omega_C,\omega_C)=0$, and by  \cite[(2.a.17)]{Gr}  it is
   $\sK_{0,1}(C,\omega_C,\omega_C)= \sK_{0,2}(C,\omega_C)$.

    Write $C= A+ \Gamma_h+\Gamma_k$:
since it is $p_a(\Gamma_i)=0$,  $p_a(
\Gamma_h+\Gamma_k )=-1$ and all these curves are not disconnecting then 
we can  consider the long exact sequence of Koszul groups (Prop.  \ref{splitting}) 
for the decomposition $C = A+\Gamma_h+\Gamma_k $  (respectively for  the decompositions
 $C- \Gamma_h = A+\Gamma_k$, $C-\Gamma_k= A+\Gamma_h$). 
 
Let  $W=\im \{ H^0(\omega_{C}) \rightarrow H^0(A, \omega_{C}) \}$.  
By Thm. \ref{teo:2} for    $i \in \{h,k\}  $  and every $q\geq 1$ it is
  $\sK_{0,q}(A +\Gamma_i, \omega_{C})=0$; consequently it is 
 $\sK_{0,q}(A,W,  \omega_{C})=0  \  \forall q\geq 1$.
 
 Therefore it is sufficient to prove that  the following sequence is exact
  \[ 
\xymatrix
{ 
\sK_{1,1} (C, \omega_{C} )   \ar[r] ^{\hspace{ - 4mm} \iota} &  \sK_{1,1} (A, W , \omega_{C} )    \ar[r] ^{ \hspace{ - 14mm} \pi}& \sK_{0,2}  (\Gamma_h+\Gamma_k, \Oh_{\Gamma_h+\Gamma_k}(-A ),\omega_{C})).
}\]
By \cite[(2.a.17)]{Gr} and adjunction  we  have 
 $$\sK_{1,1} (\Gamma_k, \Oh_{\Gamma_k}(-A -\Gamma_h), \omega_{C})  = \sK_{1,0} (\Gamma_k, \omega_{\Gamma_k} , \omega_{C}).  $$  
By 
 \cite[(2.a.17), (3.a.6)]{Gr} we get $  \sK_{1,0} (\Gamma_k, \omega_{\Gamma_k}, \omega_{C})  =0 $
 since $\Gamma_k \iso \proj^1$. 
 Similarly 
 $ \sK_{1,1} (\Gamma_h, \Oh_{\Gamma_h}(-A -\Gamma_k ), \omega_{C}) = 0 $.
Moreover, since $\Gamma_h\cap \Gamma_k= \{ \emptyset\} $ it is $\Oh_{\Gamma_k}(-A-\Gamma_h)
\iso \Oh_{\Gamma_k}(-A)$ and $\Oh_{\Gamma_h}(-A-\Gamma_k)
\iso \Oh_{\Gamma_h}(-A)$, whence we obtain 
 the splitting of the exact sequence of invertible sheaves 
$$ 0\to \Oh_{\Gamma_k}(-A-\Gamma_h) \to \Oh_{\Gamma_h +\Gamma_k} (-A) \to
 \Oh_{\Gamma_h}(-A)\to 0,$$
and for every $p,q$   the isomorphism
  $$
\sK_{p,q}  (\Gamma_h+\Gamma_k, \Oh_{\Gamma_h+\Gamma_k}(-A ),\omega_{C})) \iso  \sK_{p,q}  (\Gamma_h, \Oh_{\Gamma_h}(-A), \omega_{C})  
 \bigoplus \sK_{p,q}  (\Gamma_k, \Oh_{\Gamma_k}(-A ),\omega_{C}))
 $$
 In particular we get $\sK_{1,1}  (\Gamma_h+\Gamma_k, \Oh_{\Gamma_h+\Gamma_k}(-A ),\omega_{C})) =0$, that is $\iota$ is injective. To prove the surjectivity of $\pi$
  we consider   the following commutative diagram:
 
{ \footnotesize
  \[ 
\xymatrix@C=+3mm
{ 
\sK_{1,1} (C, \omega_{C} )  \ar@{^{(}->}[d]  \ar@{^{(}->}[r] 
 & \sK_{1,1} (A+\Gamma_h   , \omega_{C} )  \ar@{^{(}->}[d] \ar[dr]
 &  \\
 \sK_{1,1} (A  + \Gamma_k, \omega_{C} )  \ar[dr]   \ar@{^{(}->}[r] 
 & \sK_{1,1} (A, W , \omega_{C} )  \ar[d]^{\pi_{2}}  \ar[dr]^{\pi}   \ar@{->}[r] ^{ \hspace{ - 5mm} \pi_{1}} 
 & \sK_{0,2}  (\Gamma_k, \Oh_{\Gamma_k}(-A ), \omega_{C})   \ar[d] 
 \\
 & \sK_{0,2}  (\Gamma_h, \Oh_{\Gamma_h}(-A), \omega_{C})  \ar[r]
 & \sK_{0,2}  (\Gamma_h, \Oh_{\Gamma_h}(-A), \omega_{C})  \oplus  \sK_{0,2}  (\Gamma_k, \Oh_{\Gamma_k}(-A ), \omega_{C}) 
 \\
}
\]}
Now $p$ is surjective since $\sK_{0,2}( A  + \Gamma_k, \omega_{C} )=0$.  Analogously $q$ is surjective.
Moreover  we have $\pi=(\pi_1, \pi_2)$, and   $\sK_{1,1} (C, \omega_{C} ) =  \sK_{1,1} (A+\Gamma_h   , \omega_{C} ) \cap 
  \sK_{1,1} (A+\Gamma_k   , \omega_{C} ) $ (that is 
the space of  quadrics vanishing along $C$ 
   is given considering the intersection of the quadrics vanishing along $A+ \Gamma_h$, resp. along $A+\Gamma_h$).  
 Therefore $\pi$ is surjective, which implies   $\sK_{0,2} (C, \omega_{C} ) =0$.

  \hfill\break
 {\bf  {\em   C.4.}  }  
{\em  For every irreducible subcurve $\Gamma_i$  it is $\Gamma_i^2\leq 1$,  and for 
for every distinct pairs $\Gamma_i, \Gamma_j $ it is 
 $\Gamma_i \cdot \Gamma_j =1$; moreover     the curve $(\Gamma_i+\Gamma_j) \subset C$ is  always not disconnecting.  }

 \hfill\break
 {\bf  {\em   C.4.1.}  } 
Assume that $C$ contains three  components  
 $\Gamma_1, \Gamma_2 , \Gamma_3$ (possibly equal) such that 
 $\Gamma_1\cdot \Gamma_2 = \Gamma_2 \cdot \Gamma_3 =\Gamma_1\cdot \Gamma_3 =1 $ and  
 $\Gamma:= \Gamma_1+ \Gamma_2 + \Gamma_3$ is not disconnecting. 
(Notice that  if $C$ has only one irreducible component, then  we are exactly in this case since necessarily it is
$\Gamma_1^2=1$  and $\mult_C(\Gamma_1) \geq 5 $).

 In this case 
 $\Gamma$ is 2-connected  with arithmetic genus =1,   $E= C-G \neq \emptyset$ and then we  can proceed as in {\em Case A}, since it is
 $\omega_{\Gamma} \iso \Oh_{\Gamma} $.

\hfill\break
 {\bf  {\em   C.4.2.}  }  
 Assume that $C$ contains three  components  
 $\Gamma_1, \Gamma_2 , \Gamma_3$  (possibly equal)  such that 
 $\Gamma_1\cdot \Gamma_2 = \Gamma_2 \cdot \Gamma_3 =\Gamma_1\cdot \Gamma_3 =1 $    and the curve
 $\Gamma= \Gamma_1+ \Gamma_2 + \Gamma_3$ is disconnecting. 
 
 In this case 
we can write 
 $C-\Gamma_1 - \Gamma_2= E+\Gamma_3 +G'$ with $E, G'$ as in  Lemma \ref{disconnecting}, 
 that is,  we have a decomposition $C= E+ \Gamma+G'$  with $\Oh_{G'}(-E)$ NEF.   Moreover 
it is 
 $E\cdot \Gamma_3 \geq 1$, $G' \cdot \Gamma_3\geq 1$ since  $\Gamma_3\iso \proj^1$ and $\Gamma_1+\Gamma_2$ is not disconnecting, and similar
inequalities  hold for $\Gamma_1$ and $\Gamma_2$.  

Let $G= G'+ \Gamma$. Since $E$ is connected  It is  enough to prove that 
$r_1:H^0(G,{\omega_C} \otimes H^0(G,\omega_G) \to H^0(G, \omega_C \otimes \omega_G) $.

Notice that for every  $i\in \{1,2,3\}$ we have  $\deg_{\Gamma_i} \omega_{G}  \geq 0$  and 
$H^0(G,\omega_G) \onto H^0(\Gamma_i,\omega_G) $. 

Without loss of  generality we may 
assume $E\cdot \Gamma_1\geq 2$
since $E\cdot (\Gamma_1 +\Gamma_2+ \Gamma_3) \geq 4$. 
We work as in  { \em Case C.2.1}, i.e., we consider the splitting $G=\Gamma+G'$ and 
we take the sheaf $\omega_{\Gamma} (E)$.   Then we have  the exact sequence 
 $$ 0\to H^0(\Gamma_1, \omega_{\Gamma_{1}} (E)) \to
 \ H^0(\Gamma, \omega_{\Gamma} (E)) \to \ H^0(\Gamma_2 +\Gamma_3, \omega_{\Gamma_{2}+\Gamma_3} (\Gamma_1 +E)) \to 0 
 $$
and
by  the same degree arguments  adopted in  {\em Case C.2.1} it is immediately  seen that we have the surjective maps 
 $$ H^0(\Gamma_1, \omega_{\Gamma_1} (E)) \otimes H^0(\Gamma_1, \omega_G) \onto 
 H^0(\Gamma_1, \omega_{\Gamma_1} (E)\otimes \omega_G)$$
 $$H^0(\Gamma_2+\Gamma_3, \omega_{\Gamma_2 +\Gamma_3} (\Gamma_1 +E)) \otimes H^0(\Gamma_2+\Gamma_3, \omega_G) \onto 
H^0(\Gamma_2+\Gamma_3, \omega_{\Gamma_2+\Gamma_3} (\Gamma_1 +E)\otimes \omega_G)$$ 
that is, we get the surjection
$H^0(\Gamma, \omega_{\Gamma} (E)) \otimes H^0(\omega_G) \onto H^0(\Gamma, \omega_{\Gamma} (E) \otimes \omega_G)$.
Finally, as in the previous cases $H^0(G, \omega_{C}) \otimes H^0(\omega_G) \onto H^0(G, \omega_{C} \otimes \omega_G)$
since $\Oh_{G'}(-E)$ is NEF, and then we can conclude that $r_1$ is onto.

  \hfill\break
 {\bf  {\em   C.4.3.}  }  
Finally,  
 we are left with the case where $C$ has exactly two irreducible components, $\Gamma_1,\Gamma_2$ of nonpositive selfintersection:   
$C=n_1\Gamma_1 +n_2\Gamma_2$,  $\Gamma_1 \cdot \Gamma_2=1$ and $\Gamma_i^2 \leq 0$ for $i=1,2$. 

We may 
assume $\Gamma_1^2 \geq \Gamma_2^2$. Since $C$ is 4-connected  with an easy computation 
we obtain 
$\Gamma_1^2 = 0 $ and  $\mult_C(\Gamma_i)\geq 4$.  Moreover $n_2$ is even since $C$ is an even curve.

Notice that  for every subcurve $B= \alpha_1 \Gamma_1 +\alpha_2\Gamma_2 \subset C$ 
it is 
$$B\cdot (C-B)= \alpha_2(n_1 + (n_2-1)\Gamma_2^2) - \alpha_2(\alpha_2-1) \Gamma_2^2 + \alpha_1(n_2-2\alpha_2)
 \geq \alpha_2(n_1+(n_2-1)\Gamma_2^2)$$
(since we may assume  $2 \alpha_2 \leq n_2$ by the symmetry of the intersection product), which implies 
 $B\cdot (C-B)\geq 4 \alpha_2 $ since $\Gamma_2 \cdot (C-\Gamma_2)=n_1+(n_2-1)\Gamma_2^2\geq 4$. 

If $\Gamma_2^2=0$, we take $E= 2\Gamma_1+2\Gamma_2$.  Then $p_a(E)=1$ and  applying a refinement of the above formula  it is easy to see that $C-E$ is numerically connected.
In this case  we can conclude as in {\em Case A.}

If  $\Gamma_2^2<0$, let $a_1=   \lceil  \frac{n_1}{2} \rceil  $,  $a_2 = \frac{n_2}{2} $
and let 
$G:=a_1   \Gamma_1 + a_2\Gamma_2$, $E:= C-G = (n_1-a_1)\Gamma_1+a_2\Gamma_2$.

Now $E$ is numerically connected and $G$ is 2-connected.  Indeed,  let us consider a subcurve 
 $B=\alpha_1 \Gamma_1 +\alpha_2\Gamma_2  \subset G$. Since $2G\cdot B\geq G\cdot B$ 
we have
$ B\cdot (G-B) \geq \frac{1}{4} 2B\cdot (C-2B)  \geq   2$ since $2B \cdot (C-2B) \geq  8\alpha_2$ by the above formula. If $B\subset E$  it is $ B\cdot (E-B) \geq \frac{1}{4} 2B\cdot (C-\Gamma_1 -2B)  \geq   1$.  

Therefore it is enough to prove the surjection of $r_1: H^0( G, \om_{C}) \otimes  H^0( G, \om_{G}) \to 
H^0( G, \om_{G} \otimes \omega_C)$. We have  
 the following exact sequence
$$  0 \to H^0(a_1 \Gamma_1 + \Gamma_2, {\om}_{a_1\Gamma_1+\Gamma_2} (E) ) \to 
 H^0( G, \om_{C}) \to  H^0((a_2-1)\Gamma_2, \om_{C} ) \to 0$$
 and moreover $ H^0( G, \om_{G}) \onto  H^0((a_2-1)\Gamma_2, \om_{G} )$ since $a_1 \Gamma_1 + \Gamma_2$ is numerically connected.

By  \cite[Thms. 3.3]{CFHR}
 $|\omega_G|$ is a base point free system on $G$ since $G$ is 2-connected. 
 Let  
$W := \im \{ H^0(G, \omega_G) \to H^0(a_1 \Gamma_1 + \Gamma_2, \omega_G)\} $. 
Then $W$  is a base point free system  and moreover we have 
$
H^1(a_1 \Gamma_1 + \Gamma_2,  \om_{a_1\Gamma_1}(E) \otimes {\omega_G}^{-1}) \iso H^1(a_1\Gamma_1 +\Gamma_2,
 E-(a_2-1)\Gamma_2) = 0 $
because  $\Gamma_i \iso \proj^1$, $\Gamma_1^2=0$, $\Gamma_1\cdot \Gamma_2=1$ and $(E-(a_2-1)\Gamma_2) \cdot \Gamma_1 =1$, 
 $(E-(a_2-1)\Gamma_2) \cdot \Gamma_2 \geq 1 $ since $E$ is 1-connected.
Therefore by   Prop.  \ref{prop:van}  we have the surjection  
$$H^0(a_1 \Gamma_1 + \Gamma_2, {\om}_{a_1 \Gamma_1 + \Gamma_2} (E) ) \otimes W  \onto 
 H^0(a_1 \Gamma_1 + \Gamma_2 , {\om}_{a_1 \Gamma_1 + \Gamma_2} (E)  \otimes \om_{G}) $$
 
Finally $  H^0((a_2-1)\Gamma_2, \om_{G} ) \otimes  H^0((a_2-1)\Gamma_2, \om_{C} ) \onto 
 H^0((a_2-1)\Gamma_2, \om_{G}\otimes\om_C  )$ follows from 
({\em i})  of Prop. \ref{prop:van}  taking $\sH =\omega_C, \sF=\omega_G$ if $\Oh_{\Gamma_2}(E)$ is NEF,  
or $\sF=\omega_C, \sH=\omega_G$ if $\Oh_{\Gamma_2}(-E)$ is NEF. 

	\hfill{   \em Q.E.D.  for  Theorem \ref{teo:3}}

\section{ On the canonical ring of regular surfaces}
In this section we prove Theorem \ref{teo:bicanonical}.
The arguments we adopt are  very classical and based on the ideas 
 developed in  \cite{ci}.  Essentially we simply  restrict to a curve in the canonical system $|K_S|$. 
The only novelty is that now we do not
make any  requests on such a curve
 (i.e. we allow the curve $ C\in |K_S|$ to be singular and with many components)  since we can apply 
 Thm.  \ref{teo:3}.

\hfill\break
	{\bf Proof of  Theorem \ref{teo:bicanonical}.} 
	
	By assumption there exists a   $3 $-connected not 
	honestly hyperelliptic curve 
	 $C =\sum_{i=1}^{s} n_{i} \Gamma_{i} \in |K_{S}|$. 
	 Let $s \in H^0(S, K_S ) $ be the corresponding  section, so that $C$ 
	 is   defined
 by $(s)=0$.  

	By adjunction we have 
	$(K_{S}^{\otimes 2})_{|C}= 
(K_{S}+C)_{|C} \iso 
\omega_{C}$, that is, $C$  is an  even curve; in particular 
it is   $4$-connected. 
Thus we can apply  Theorem \ref{teo:3}, obtaining  the surjection 
$$ (H^{0}(C,K_{S}^{\otimes 2}) )^{\otimes k} \onto  H^{0}(C,K_{S}^{\otimes 2k})  \ \ \ \ \forall k \in \Na. $$ 
Now let  us consider 
the usual maps given by multiplication of sections
\begin{eqnarray*}
A_{l,m}: H^0(S, K_S^{\otimes l}) \otimes H^0(S, K_S^{\otimes m}) \rightarrow H^0(S, K_S^{\otimes (l+m)})\\  
a_{l,m}: H^0(C, K_S^{\otimes l}) \otimes H^0(C, K_S^{\otimes m}) \rightarrow H^0(C, K_S^{\otimes (l+m)})
\end{eqnarray*}
and consider
 the following  commutative  diagram 
   \[ 
\xymatrix{ 0
 \ar[d] & 0 \ar[d] \\ 
H^0(S, K_S^{\otimes (k -1)} ) \ar[d]^{C_{k}} \ar[r]^{\iso} & H^0(S, K_S^{\otimes (k -1)})  \ar[d]^{c_{k}} \\
{\displaystyle
\bigoplus_{{\tiny{\begin{array}{l} 
l+m=k \\
[-2 pt] 
0<l\leq m 
\end{array}}}}}
\left[ H^0(S, K_S^{\otimes l}) \otimes H^0(S, K_S^{\otimes m})\right]  \ar[d]^{R_{k}}  \ar[r]^{\hspace{ 2cm}{\overline{{\rho}}_{k}}} &
H^0(S, K_S^{\otimes k})\ar[d]^{r_{k}} \\ 
{\displaystyle
\bigoplus_{{\tiny{\begin{array}{l} 
l+m=k \\
[-2 pt] 
0<l\leq m 
\end{array}}}}}
\left[ H^0(C, K_S^{\otimes l}) \otimes H^0(C, K_S^{\otimes m})\right]  \ar[d] \ar[r]^{\hspace{ 2cm}{{\rho}_{k}}}&
H^0(C, K_S^{\otimes k})\ar[d] \\ 
0 & 0 
}
\]
Here  the left hand column is a  complex,  while the right hand column is exact. Moreover 
\begin{enumerate}
\item[$\bullet$] 
$C_k$ is given by tensor product with $s$ while  $c_k$ is given by product with $s$ 
\item[$\bullet$]  $R_k =
 {\displaystyle
\bigoplus_{{\tiny{\begin{array}{l} 
l+m=k \\
[-2 pt] 
0<l\leq m 
\end{array}}}}} r_l \otimes r_m $ \ \  (where $r_l, r_m$ are the usual restrictions)
\item[$\bullet$]
 ${\overline{\rho}}_k =  {\displaystyle\bigoplus_{{\tiny{\begin{array}{l} 
l+m=k \\
[-2 pt] 
0<l\leq m 
\end{array}}}}} A_{l,m}  $  \ \ and \ \  
${\rho}_k =  {\displaystyle \bigoplus_{{\tiny{\begin{array}{l} 
l+m=k \\
[-2 pt] 
0<l\leq m 
\end{array}}}}} a_{l,m}  $\ \  
\end{enumerate}
Note that it is $ \coker ({\rho}_{k}) \iso \coker ({\overline{{\rho}}_{k}})$ for every $k \in \Na$. 

Now, for $S$ of general type,  if $p_g\geq 1$ and $q=0$ by \cite[Thm. 3.4 ]{ci} 
${\overline{{\rho}}_{k}}$ is surjective for every
$k \geq 5$ except  the case  $p_g=2, \ K^2=1$, which is not our case since 
otherwise $C$ would be a curve of genus $2$ contradicting our assumptions.
For $k=4$  the map  $a_{2,2}$  is surjective by  Thm. \ref{teo:3}. Whence   $\rho_4$ and $\overline{\rho}_4$
are surjective, and this proves the theorem.

\hfill{   \em Q.E.D.  for  Theorem \ref{teo:bicanonical}.}

\noindent {\bf Author's address:}

\bigskip

\noindent Marco Franciosi \\ Dipartimento di Matemativa Applicata ``U.
Dini", Universit\`a di Pisa \\
    via Buonarroti 1C, I-56127, Pisa, Italy \\  e-mail:
franciosi@dma.unipi.it

\end{document}